\input ./style/arxiv-general.cfg
\documentclass[MSNbibl,number,citesort,seceqn,dvips]{arxbj}
\makeatletter
   \@ifpackageloaded{graphicx}{}{\usepackage{graphicx}}
\makeatother
\usepackage{multirow}

\volume{22}
\issue{4}
\pubyear{2016}
\firstpage{2521}
\lastpage{2547}
\doi{10.3150/15-BEJ736}
\docsubty{FLA}

\makeatletter
\renewcommand{\mid}{|}
\newcommand{\rrvert}{\vert}
\newcommand{\rrVert}{\Vert}
\newcommand{\llvert}{\vert}
\newcommand{\llVert}{\Vert}
\newtheorem{teo}{Theorem}[section]
\newtheorem{lem}{Lemma}[section]
\newremark{Remark}{Remark}[section]
\newtheorem{cor}{Corollary}[section]
\newremark{ex}{Example}[section]
\newproclaim{d1}{Definition}[section]
\newtheorem{pro}{Proposition}[section]
\newcommand{\kk}{\kappa}
\newcommand{\ttt}{\theta}
\newcommand{\aaa}{\alpha}
\newcommand{\bb}{\beta}
\newcommand{\bex}{
\begin{ex}}
\newcommand{\eex}{\end{ex}}

\newcommand{\bp}{
\begin{pro}}
\newcommand{\ep}{\end{pro}}
\newcommand{\be}{
\begin{equation}
} \newcommand{\ee} {
\end{equation}}
\newcommand{\bl}{\begin{Lemma}}
\newcommand{\el}{\end{Lemma}}
\newcommand{\bn}{\begin{notation}}
\newcommand{\en}{\end{notation}}
\newcommand{\bcon}{\begin{construction}}
\newcommand{\econ}{\end{construction}}
\newcommand{\bt}{\begin{theorem}}
\newcommand{\et}{\end{theorem}}
\newcommand{\br}{
\begin{Remark}}
\newcommand{\er}{\end{Remark}}
\newcommand{\bc}{\begin{corollary}}
\newcommand{\ec}{\end{corollary}}

\newcommand{\normt}[1]{\llVert #1\rrVert}
\makeatother

\begin{document}
\begin{frontmatter}

\title{Methods for improving estimators of  truncated circular parameters}
\runtitle{Truncated circular parameters}

\begin{aug}
\author[A]{\inits{}\fnms{}\hspace*{-1pt}\snm{Kanika}\thanksref{e1}\ead[label=e1,mark]{10ma90r04@iitkgp.ac.in}}
\and
\author[A]{\inits{S.}\fnms{Somesh}~\snm{Kumar}\corref{}\thanksref{e2}\ead[label=e2,mark]{smsh@maths.iitkgp.ernet.in}}
\address[A]{Department of Mathematics, Indian Institute of Technology Kharagpur,  Kharagpur-721302, West Bengal, India. \printead{e1,e2}}
\end{aug}

%
\received{\smonth{1} \syear{2015}}

\begin{abstract}
In decision theoretic estimation of parameters in  Euclidean space $\mathbb{R}^p$, the  action
space is chosen to be the convex closure of the estimand space. In this paper, the
concept has been extended to  the estimation of circular parameters of distributions
having support as a circle, torus or cylinder.    As directional distributions are of
curved nature, existing methods for distributions with parameters taking values in
$\mathbb{R}^p$ are not immediately  applicable here.  A circle is the simplest one-dimensional
Riemannian manifold. We employ concepts of convexity, projection, etc., on manifolds to
develop sufficient conditions for inadmissibility of estimators for circular parameters.
Further invariance under a compact group of transformations is introduced in the estimation
problem   and a complete class theorem  for equivariant estimators is derived.
This extends the results of Moors [\textit{J.~Amer. Statist. Assoc.} \textbf{76} (1981) 910--915]
on $\mathbb{R}^p$ to circles.
The findings are of special interest to the case when a circular parameter is truncated.
The results are implemented to a wide range of directional distributions to obtain  improved  estimators of circular parameters.
\end{abstract}

\end{frontmatter}

\section{Introduction}
Problems of  estimation when  the parameter space is restricted  are  encountered often  in practice.
These restrictions arise due to prior information on parameters and they
can be in the form of bounds on the range  or  equality/inequality constraints  of several  parameters. For recent developments and discussions on various aspects of  estimation procedures in restricted parameter space problems, one may refer to
\cite{kub04,kub05,marstr04,marstr12,ve06} and references therein. Frequently in
practical  applications, we assume the random observations taking values in Euclidean
spaces. However,  it sometimes may be more useful to represent them on circles/spheres/cylinders. In such cases, we employ directional distributions.  For instance, mortality data due to a specific   disease may be better represented as circular data to study the seasonal pattern of the   disease.  There are numerous situations in biological,  meteorological,  astronomical  applications, where    directional data (circular/axial/spherical)   arises  \cite{b81,f93,mbok}. However,   little  attention  has been paid to  problems of estimating directional parameters under constraints.
Rueda, Fern{\'a}ndez and Peddada \cite{rud09} considered the estimation of the circular  parameters under order restrictions. There are situations when the  parameter may lie  on an  arc of the circle.  For example,  the peak of mortality rates due to respiratory diseases occurs  during November  to February.

One  major consequence of placing restrictions on the natural parameter space  is that  estimators derived using standard concepts of maximum likelihood, minimaxity, invariance, etc., become inadmissible. However, existing methods  developed for Euclidean spaces  $\mathbb{R}^{p}$ are not  directly   applicable to directions which are represented to lie  on unit hypersphere with the center at origin $\mathbb{S}_{p-1}=\{\mathbf{x} \in\mathbb{R}^{p}: \normt{\mathbf{x}}=(\mathbf{x}^{T}\mathbf{x})^{1/2}=1\}$. Topological properties of $\mathbb{S}_{p-1}$ depend on the  differential geometry of  embedding $\mathbb{S}_{p-1}$ in  $\mathbb{R}^p$ as  $\mathbf{x}\rightarrow{\mathbf{x}}/{\normt{\mathbf{x}}}$ with $\mathbf{x}\in\mathbb{R}^p$. We need to suitably modify techniques available for $\mathbb{R}^{p}$  to improve  standard  estimators of directions. In this paper, we consider the case of $p=2$ and denote by $\mathbb{S}$ a unit  circle. The elements of  $\mathbb{S}$ can be specified by corresponding angles with respect to an arbitrary choice of  zero direction and  orientation.  Let us define by  $\mathbb{T}=[0,2\pi)$ the space of  amplitudes  of the unit vectors in  $\mathbb{S}$. The\vspace*{1pt} point on  $\mathbb{S}$ corresponding to an angle $\aaa\in\mathbb{T}$ is $(\cos \alpha, \sin \alpha)^{T}$ and the angle corresponding to  a point  $\mathbf{x}=(x_1,x_2)^T\in\mathbb{S}$ is $\mathsf{atan} ({x_2}/{x_1} )$, where  the function $\mathsf{atan}(\cdot )$ is given in (\ref{atan}). Note that $\mathbb{S}$ and $\mathbb{T}$ are isomorphic.

Let  random variable $\mathbf{Z}\in\mathfrak{Z}$  have  an unknown probability distribution $\mathsf{P}_{\bolds{\nu}}, \bolds{\nu}\in \Omega$. Let the family $\{\mathsf{P}_{\bolds{\nu}}, \bolds{\nu}\in \Omega\}$ be dominated by a measure $\eta$. The support $\mathfrak{Z}$ may be circle $\mathbb{T}$, torus $\mathbb{T}^k$ or cylinder $\mathbb{R}\times\mathbb{T}$, etc.; however, the  estimand $h(\bolds{\nu})$ is circular, where $h(\cdot )$ is a measurable function from  $\Omega$ into $\mathbb{T}$.    The problem  of estimating $h(\bolds{\nu})$ is considered   under a circular loss function:
\begin{eqnarray}
\label{ch5_loss} \mathsf{L}(\bolds{\nu},{\delta})=1-\cos \bigl(h(\bolds{ \nu}) -
\delta \bigr).
\end{eqnarray}
In the case when the parameter space is a subspace of  $\mathbb{R}^{p}$  and the loss function is an  increasing function of  Euclidean distance, the action space is chosen  as a  convex closure of the range of estimand.  Unlike this well-known result,  it is demonstrated in Section~\ref{SE3}   that analogous result does not necessarily hold for  circular parameter $h(\bolds{\nu})$ under the loss $\mathsf{L}$.

 One of the major contributions was of Moors \cite{moors81} (see \cite{moors85}, Chapter~3, also) to the estimation problem of truncated parameters  of a  unknown family of distributions on  $\mathbb{R}^{p}$  dominated by a  $\sigma$-finite measure. Estimators (except the  constant ones)  taking value near the boundaries of   action space with the  positive probability  turn out to be inadmissible with respect to the  squared loss function under certain conditions on the transformation group. He considered invariance  under a finite group with measure preserving elements such that  induced transformations of the action space satisfy the  linearity  property and group of these transformations is commutative. Under this scenario, he  constructed a subspace of the original action space  and proved that any invariant estimator taking values outside this new action space with the  positive probability is inadmissible  and dominated by its projection on the new action space.  Later,
Moors and van Houwelingen
\cite{moors93} relaxed  conditions of measure preserving and commutativity. Without dropping these   conditions,
Kumar and Sharma \cite{ks92}  generalized the result of Moors \cite{moors81} to a locally compact group  such that  induced transformations  on the action space are affine (stronger condition than linear property) and  loss function is an  increasing function of Euclidean distance.  Along with these ideas, an analogous theory for  circular parameter $h(\bolds{\nu})$ is developed in Section~\ref{SE4}.

 The outline of the paper is as follows.  The concepts of distance formulae,  convexity, closure of a set and  projections play a prominent  role here. Section~\ref{SE2}  provides the mathematical background of these concepts for  $\mathbb{T}$.  In Section~\ref{SE3}, we consider the estimation of  circular parameter $h(\bolds{\nu})$ when it is restricted to lie on an  arc of circle and estimation space $\mathcal{A}$ is chosen as the convex closure of $h(\Omega)$. A complete class result for this estimation problem is obtained under certain conditions. Then the  result is illustrated for several  directional distributions.   In Section~\ref{SE4}, we introduce invariance under a compact group $\mathcal{G}$ in the  estimation problem such that induced  transformations on   $\mathcal{A}$ satisfy  the circular  property.  Sufficient condition for inadmissibility of an $\mathcal{G}$-equivariant estimator is obtained. Applications of this  result are demonstrated for both unrestricted or restricted estimation problems. For restricted estimation problems, improved estimators obtained in Section~\ref{SE3} are further improved using the result of  Section~\ref{SE4}.

\section{Definitions and preliminary results}
\label{SE2}
Before we embark on estimation problem, we introduce some preliminary results in this section. For a subset  $A\subset\mathbb{T}$,  Lebesgue measure,  interior,  convex hull,  convex closure and boundary  of   $A$ are denoted by $l(A)$, $\mathsf{int}(A)$,   $\mathsf{conv}(A)$, $\mathsf{cc}(A)$ and $\mathsf{bd}(A)$, respectively.

\subsection{Convexity}
It is   more convenient to deal with polar coordinates than Cartesian coordinates  when observations lie   on a unit circle. We summarize   the   concept of convexity  for the circle  $\mathbb{S}$ and  then adopt it  for the space $\mathbb{T}$.

 A geodesic (\cite{cu94}, page~15), on Riemannian manifold generalizes the line in Euclidean space.    In the context of Riemannian manifold $\mathbb{M}$ equipped with Riemannian metric,  a subset $A\subset\mathbb{M}$ is convex if minimal geodesic with the end points in $A$ belongs to $A$.  In the case of $\mathbb{S}$, great circles are geodesics.    Minimal geodesic between  any two points on $\mathbb{S}$ is unique unless points are antipodal (diametrically opposite). Some concepts of convexity on $\mathbb{S}$ were introduced in \cite{d63}, Section~9.1. Here, we  use  convexity and strong convexity as given below.

\begin{d1}[(Convex)]
 A set $A\subset\mathbb{S}$ is convex if for any two points in $A$, there exists a minor arc  of  great circle lying  entirely in $A$ joining them.
\end{d1}

By convention, this  definition allows antipodal points to lie  in convex sets.  Every segment of a semicircle is a convex subset of $\mathbb{S}$.

\begin{d1}[(Strongly convex)]
A subset $A$ of $\mathbb{S}$ is strongly convex if  $A$ is convex and does not contain antipodal points.
\end{d1}
In the case of  $\mathbb{R}^{p}$, convex hull of any subset is the collection  of all possible weighted arithmetic mean of   points in that subset. Analogously  for a  set $A\subset\mathbb{S}$,  the convex hull of $A$  is the smallest convex set (not necessarily strong convex)  containing $A$, that is,  it consists of
\begin{eqnarray}
\label{convex_comb} {(w_1\mathbf{x}_1+ \cdots+w_n
\mathbf{x}_n)}/{\normt{w_1 \mathbf{x}_1+
\cdots+w_n\mathbf{x}_n}}
\end{eqnarray}
for all nonnegative weights $w_1,\ldots,w_n$ such that $\sum_{i=1}^{n}w_i=1$ and  for all $\mathbf{x}_1,\ldots,\mathbf{x}_n\in A$ provided   the norm in the denominator is nonzero. The convex hull of  two antipodal points of $\mathbb{S}$  does not exist (\cite{b01}, Section~2.3). Polar form of (\ref{convex_comb}) is discussed in Section~\ref{CE}.

Extension of a   fundamental theorem of Carath\'eodory for $\mathbb{R}^{2}$  to $\mathbb{S}$   can be stated as below (see~\cite{b01}).
\begin{lem}\label{convex_4}
Each point in the convex hull of a set $A\subset\mathbb{S}$ can be expressed as normalized weighted arithmetic mean of at most $2$ points of  $A$.
\end{lem}

To extend definitions of convex and strong convex sets to the space  $\mathbb{T}$,  we need   certain  subsets of $\mathbb{T}$.  For  $\aaa, \bb\in\mathbb{T}$,  definitions for sets  of type $I$, $J$, $K$  and $K_1$ are stated only for intervals of form $(\aaa,\bb)$.  They  may be extended to intervals of other
forms $(\aaa,\bb],  [\aaa,\bb)$ and  $[\aaa,\bb]$.  For  $\aaa\leq\bb$, let
\begin{eqnarray*}
I(\aaa,\bb)&=& (\aaa,\bb),
\\
J(\aaa,\bb)&=& [0,\aaa)\cup(\bb,2\pi),
\\
K(\aaa,\bb)&=& \cases{ I(\aaa,\bb), &\quad if $0\leq(\bb-\aaa)<\pi$;
\cr
I(\aaa,
\bb)\mbox{ or }J(\aaa,\bb), &\quad if $(\bb-\aaa)=\pi$;
\cr
J(\aaa,\bb), &\quad if
$\pi<(\bb-\aaa)< 2\pi$.}
\end{eqnarray*}
  Sets  $I[\aaa,\bb]$ and $J[\aaa,\bb]$ contain all the angles corresponding to an arc joining two  points $(\cos\aaa,\sin\aaa)^T$ and $(\cos\bb,\sin\bb)^T$  in the positive and the negative directions, respectively. Moreover,   $K[\aaa,\bb]$  contains angles corresponding  to the minor arc  joining them.
\br
Although in  definition of $I$-type set, $\bb$ is not allowed to  take value $2\pi$,  intervals  $(\aaa,2\pi)$ and $[\aaa,2\pi)$ can be expressed as $J(0,\aaa)$ and $J(0,\aaa]$, respectively, for all $\aaa\in\mathbb{T}$. For any $\aaa\in\mathbb{T}$, $J(\aaa,\aaa)=\mathbb{T}-\{\aaa\}$ and  $J(\aaa,\aaa]=J[\aaa,\aaa)=J[\aaa,\aaa]=\mathbb{T}$.
\er
For defining sets of type $I$, $J$ and $K$, we have taken $\aaa\leq\bb$. The following definition for $K_1$-type subsets  of $\mathbb{T}$ does not have this restriction.
\begin{eqnarray*}
K_1(\aaa,\bb)&=&\cases{ K(\aaa,\bb), &\quad if $\aaa\leq\bb$,
\cr
K(
\bb,\aaa), &\quad if $\aaa>\bb$.}
\end{eqnarray*}
Note that   set $K(\aaa,\bb)$  is isomorphic to $I(0,\gamma)$ with $\gamma\in[0,\pi]$. For  $0\leq(\bb-\aaa)\leq\pi$ and $\pi<(\bb-\aaa)<2\pi$, $K(\aaa,\bb)$ can be  transformed  to  $I(0,\bb-\aaa)$  and $I(0,2\pi-\bb+\aaa)$  using  rotation by   angles $2\pi -\aaa$ and $2\pi -\bb$,  respectively. Extending this argument, we have the following result.
\begin{lem}\label{rot}
Sets $K_1(\aaa,\bb), K_1(\aaa,\bb], K_1[\aaa,\bb)$  and $K_1[\aaa,\bb]$  with $\aaa,\bb\in\mathbb{T}$ are isomorphic  to $I(0,\gamma), I(0,\gamma], I[0,\gamma)$  and $I[0,\gamma]$, respectively,  with  $0\leq\gamma\leq\pi$.
\end{lem}
For studying topological properties, we consider the  metric space $(\mathbb{T},d)$   with the following definition of metric:
  \begin{eqnarray*}
d(\alpha, \beta)&=&1-\cos(\alpha - \beta),\qquad \aaa, \bb\in\mathbb{T}.
\end{eqnarray*}
Here, $\sqrt{2d}$ simply returns  lengths of chord   between points $(\cos\aaa,\sin\aaa)^T$ and $(\cos\bb,\sin\bb)^T$, respectively.  Consider the following  classes of subsets of $\mathbb{T}$:
\begin{eqnarray*}
\mathfrak{C}_1&=& \{\mathbb{T} \}\cup \bigl\{ K_1(
\aaa,\bb), K_1(\aaa,\bb], K_1[\aaa,\bb),
K_1[\aaa,\bb]: \aaa,\bb, \in\mathbb{T} \bigr\},
\\
\mathfrak{C}_2&=& \bigl\{K_1(\aaa,\bb),
K_1(\aaa,\bb], K_1[\aaa,\bb): \aaa,\bb
\in\mathbb{T} \bigr\} \cup \bigl\{K_1[\aaa_1,
\bb_1]: \aaa_1,\bb_1 \in\mathbb{T}\mbox{
and }\llvert \aaa_1-\bb_1\rrvert \neq\pi \bigr\},
\\
\mathfrak{C}_3&=&\mathfrak{C}_1-\mathfrak{C}_2=
\{\mathbb{T} \}\cup \bigl\{K_1[\aaa,\bb]: \aaa,\bb,\in\mathbb{T}
\mbox{ and } \llvert \aaa-\bb\rrvert =\pi \bigr\},
\\
\mathfrak{C}_4&=& \{\varnothing,\mathbb{T} \}\cup \bigl\{ I[\aaa,
\bb],J[\aaa,\bb]:\aaa,\bb,\in\mathbb{T}\mbox{ and }\aaa\leq \bb \bigr\},
\\
\mathfrak{C}_5&=&\mathfrak{C}_1\cap
\mathfrak{C}_4= \{\varnothing,\mathbb{T} \}\cup \bigl
\{K_1[\aaa,\bb]: \aaa,\bb \in\mathbb{T}\mbox{ and }\aaa\leq \bb \bigr
\}.
\end{eqnarray*}

\br
Classes $\mathfrak{C}_1$, $\mathfrak{C}_2$ and $\mathfrak{C}_5$ consist of all convex, strongly convex and  closed   convex subsets of $\mathbb{T}$, respectively. Elements of $\mathfrak{C}_3$ and $\mathfrak{C}_4$  are  closed subsets of $\mathbb{T}$. Moreover,   sets belonging to  $\mathfrak{C}_3$ (except  $\mathbb{T}$) and  $\mathfrak{C}_4$ are  corresponding to any minor arc     and  any arc on the unit circle, respectively.
\er

\subsection{Circular mean direction}\label{CE}
Since the arithmetic mean is not a suitable measure of central tendency for the angular data,   the circular mean direction is used (\cite{js01}, page~13). The weighted circular mean direction of the observations $\phi_{1},\ldots,\phi_{n}$ belonging to $\mathbb{T}$ with  weights $w_{1},\ldots,w_{n}$, $w_{i}\geq 0$ $(i=1,\ldots,n)$ such that $\sum_{i=1}^{n}w_{i}=1$ is defined  as
\begin{eqnarray*}
\bar{\phi}_w=\mathsf{atan} \biggl(\frac{\sum_{i=1}^{n}w_{i} \sin\phi_{i}}{\sum_{i=1}^{n}w_{i} \cos\phi_{i}} \biggr),
\end{eqnarray*}
with the following definition of $\mathsf{atan}(\cdot )$
\begin{eqnarray}
\label{atan} \mathsf{atan} \biggl(\frac{s}{c} \biggr) = \cases{
\tan^{-1}(s/c), &\quad if $c>0$, $s\geq0$;
\cr
\pi+\tan^{-1}(s/c),&
\quad if $c<0$;
\cr
{\pi}/{2},&\quad if $c=0$, $s>0$;
\cr
2\pi+
\tan^{-1}(s/c),& \quad if $c\geq0$, $s<0$;
\cr
\mbox{not defined},&\quad
if $s=0$, $c=0$,}
\end{eqnarray}
where $\tan^{-1}(\cdot )$ is the standard  inverse tangent  function taking values in
$ (-{\pi}/{2},{\pi}/{2} )$. The definition of $\mathsf{atan}(\cdot )$   function ensures  the following property.
\begin{lem}\label{propp1}
For all $\phi\in\mathbb{T}$ and $a\in\mathbb{R}$, we have
\begin{eqnarray*}
\mathsf{atan} \biggl(\frac{\tan\phi+a}{1-a\tan\phi} \biggr)=\bigl\{\phi+\tan^{-1}(a)
\bigr\}\bmod(2\pi).
\end{eqnarray*}
\end{lem}

Note that $\bar{\phi}_w$ is polar form of (\ref{convex_comb})   if $\phi_i$ is  corresponding angle to  $\mathbf{x}_i\in\mathbb{S}$ for all $i=1,\ldots,n$.  The following proposition proves that convex combination (weighted circular mean direction) of  finite  collection of the points in convex subset of $\mathbb{T}$   is again in that subset.

\begin{pro}\label{convex_pr}
Let   $A\in\mathfrak{C}_1$  (convex), $\phi_{1},\ldots,\phi_{n} \in A$ and  $w_{1},\ldots,w_{n}$ be nonnegative weights with $\sum_{i=1}^{n}w_{i}=1$. Then weighted circular  mean direction $\bar{\phi}_w$ of these observations belongs to $A$, if it is defined.
 \end{pro}

Circular mean direction of  a circular random variable  $\ttt$ is defined as
\begin{eqnarray*}
\mathsf{CE}(\ttt)=\mathsf{atan} ({\mathsf{E}\sin\ttt}/{\mathsf{E} \cos\ttt} ).
\end{eqnarray*}
An elementary result given in  \cite{f67}, page~74,  states  that  if random variable $\mathbf{X}$ lies in convex subset of $\mathbb{R}^p$ with probability one,  $\mathsf{E}(\mathbf{X})$ lies in the same  subset. An analogous  result for  $\mathbb{T}$  is given below.

\begin{pro}\label{jen}
If $A\in\mathfrak{C}_1$ (convex) and $\ttt$ is a random angle such that $\Pr(\ttt\in A)=1$,   the mean direction $\mathsf{CE}(\ttt)\in A$ if  $\mathsf{CE}(\ttt)$ exists. Furthermore, if $A\in\mathfrak{C}_2$ (strongly convex),  $\mathsf{CE}(\ttt)$ always exists. If $A\in\mathfrak{C}_3$ (convex but not strongly), $\mathsf{CE}(\ttt)$ does not necessarily   exist.
\end{pro}

\subsection{Projection}
The concept of projection in $\mathbb{R}^p$ is adopted to define projections of angles in $\mathbb{T}$.
\begin{d1}
The projection of an angle $\phi\in \mathbb{T}$ on a nonempty set $A\in \mathfrak{C}_4$ (closed) is defined to be the unique point $\phi_0\in A$ such that
\begin{eqnarray*}
d(\phi,\phi_0)=\inf_{\psi\in A} d(\phi,
\psi).
\end{eqnarray*}
\end{d1}
The case when $A=\mathbb{T}$ is trivial. For $\aaa,\bb\in\mathbb{T}$, let $\gamma={(\aaa+\bb)}/{2}$. If  $A$ is of from $I[\aaa,\bb]$, $\phi_0$  is given by
\begin{eqnarray*}
\phi_{0}= \cases{ \phi, &\quad if $\phi \in I[\aaa,\bb]$;
\cr
\aaa, &
\quad if $\phi \in K_1 (\aaa,\pi+\gamma ]$;
\cr
\bb, &\quad if $\phi
\in K_1 (\bb,\pi+\gamma )$;} \quad\mbox{or}\quad \cases{ \phi, &\quad
if $\phi \in I[\aaa,\bb]$;
\cr
\aaa, &\quad if $\phi \in K_1 (\aaa,
\pi+\gamma )$;
\cr
\bb, & \quad if $\phi \in K_1 (\bb,\pi+\gamma ]$.}
\end{eqnarray*}
Note that the two definitions are equivalent except when $\phi=\pi+\gamma$.  For $\phi=\pi+\gamma$,  first and second ones  yield $\phi_0=\aaa$ and  $\phi_0=\bb$, respectively. This is so because $d (\aaa,\pi+\gamma )=d (\bb,\pi+\gamma )$. If $A$ is the form of $J[\aaa,\bb]$,  $\phi_0$  is given by
\begin{eqnarray*}
\phi_{0}= \cases{ \phi, &\quad if $\phi \in J[\aaa,\bb]$;
\cr
\aaa, &
\quad if $\phi \in I \bigl(\aaa,\gamma ]$;
\cr
\bb, &\quad if $\phi \in I (\gamma,
\bb )$;} \quad\mbox{or}\quad \cases{ \phi, &\quad if $\phi \in J[\aaa,\bb]$;
\cr
\aaa, &\quad if $\phi \in I (\aaa,\gamma )$;
\cr
\bb, &\quad if $\phi \in I [
\gamma,\bb )$.}
\end{eqnarray*}
Once again the two definitions are equivalent except when $\phi= \gamma$.

Let $A$ be  a closed convex subset of  $\mathbb{R}^p$ and $\mathbf{x}\in\mathbb{R}^p$. The projection $\mathbf{x}_0$ of  $\mathbf{x}\notin  A$ on $A$  satisfies
\begin{eqnarray*}
\normt{\mathbf{x}_0-\mathbf{y}}<\normt {\mathbf{x}-\mathbf{y}}\qquad
\mbox{for all }\mathbf{y}\in A.
\end{eqnarray*}
An analogous statement holds  only  for specific closed convex subsets of $\mathbb{T}$. The following result can be easily  proved using  geometrical arguments.
\begin{lem}\label{aaaa}
   Let $\phi_0$ be the projection of an angle $\phi\notin A$  on a set $A\in \mathfrak{C}_5$ (closed convex). The inequality
\begin{eqnarray*}
d(\phi_0,\psi)< d(\phi,\psi)\qquad\mbox{for all }\psi\in A
\end{eqnarray*}
holds iff $A=I[\aaa,\bb]$ with $\bb<\aaa+(2/3)\pi$ or $A=J[\aaa,\bb]$ with $\bb>\aaa+(4/3)\pi$, that is, $l(A)<(2/3)\pi$. Moreover, if $l(A)=(2/3)\pi$ and $\mathsf{bd}(A)=\{b_1,b_2\}$, the above inequality remains strict for $\psi\in \mathsf{int}(A)$ and at least one of $\psi=b_i$ $ (i=1,2)$.
\end{lem}
For remaining sets in $\mathfrak{C}_5$, the above result holds for the expected values under certain conditions on the distribution of random variable $\ttt$ (Lemma~\ref{proj1}).

If the distribution of a circular  random variable ${\ttt}$ is  symmetric about $\psi$, the density of ${\ttt}$ with respect to any measure (measure is always finite as $\mathbb{T}$ is a compact space)  would be a function of $\cos(\ttt-\psi)$.  Further, if this distribution  is unimodal,   mean direction and mode coincide. Let us denote by    $f(\ttt\mid \psi)=f(\cos(\ttt-\psi))$  the density of a symmetric unimodal   distribution with mode $\psi$. Now consider the mixture of two unimodals $f(\cdot \mid \psi)$ and $f(\cdot \mid \psi+\pi)$. This mixture would necessarily not be bimodal.  Let $\ttt$ have a mixture distribution with probability density $\varepsilon f(\ttt\mid \psi)+(1-\varepsilon)f(\ttt\mid \psi+\pi)$, $\varepsilon\in[0,1]$. For this mixture distribution, define
\begin{eqnarray}
\label{zee} \zeta(t)=f'(t)/f'(-t),\qquad t\in[-1,1]
\end{eqnarray}
 with $t=\cos(\ttt-\psi)$. Maxima and minima of $\zeta(t)$ are denoted by   $\zeta_{\max}$  and $\zeta_{\min}$, respectively. Distribution of  $\ttt$ would be unimodal with modes $\psi$ and $\psi+\pi$ for $\varepsilon\in [\{1+\zeta_{\min}\}^{-1},1]$ and $\varepsilon\in [0,\{1+\zeta_{\max}\}^{-1}]$,  respectively. For remaining values of $\varepsilon$, it would be bimodal.

\begin{lem}\label{proj1}
Suppose that  ${\ttt}$ is a  continuous  circular  random variable whose distribution is symmetric about one of  its mode $\psi$, where  $\psi$ belongs to $A\in \mathfrak{C}_5$ (closed convex) with  $l(A)\in(2\pi/3,\pi]$  such that $\Pr (\ttt\notin A )>0$. Let  $\ttt_0$ be the projection of  ${\ttt}$  on $A$.  Then
\begin{eqnarray*}
\mathsf{E}_\psi^{\ttt} \bigl\{d (\ttt_0,\psi )
\bigr\}<\mathsf{E}_\psi^{\ttt} \bigl\{d (\ttt,\psi ) \bigr\}
\end{eqnarray*}
if distribution of ${\ttt}$  satisfies one of the following conditions:
\begin{longlist}[(C2)]
  \item[(C1)] distribution is unimodal with mode $\psi$;
  \item[(C2)] distribution is mixture with probability density  $\varepsilon f(\cdot \mid \psi)+(1-\varepsilon)f(\cdot \mid \psi+\pi)$,
  where $\varepsilon\geq1/2$ and  $\zeta(\cdot )$ defined in (\ref{zee}) is an increasing function.
  \end{longlist}
\end{lem}
It may be noted that the condition \textup{(C2)}  implies  (C1) for $\varepsilon\in[\{1+\zeta(-1)\}^{-1},1]$.
\br\label{zee1}
Convexity of density function $f(t)$  in $t\in[-1,1]$ yields increasing nature of the function $\zeta(t)$.
\er

\section{Improving estimators in  restricted parameter spaces}\label{SE3}

In Euclidean spaces, the action space  is chosen as a  convex closure of the estimand  space since estimators outside this space  with the  positive probability are  dominated by their  projections on it.  An analogous result stated below  for estimating the circular parameter ${h}(\bolds{\nu})$ is an immediate consequence of   Lemmas~\ref{aaaa} and~\ref{proj1}.

\begin{teo}\label{action_space}
Let estimand be  ${h}(\bolds{\nu})\in\Omega_1= {h}(\Omega)\subset\mathbb{T}$ and  the loss function be $\mathsf{L}$ defined in (\ref{ch5_loss}). Denote the estimation space by $\mathcal{A}=\mathsf{cc}(\Omega_1)$. Any estimator ${\delta}(\mathbf{Z})$ satisfying ${\Pr}_{\bolds{\nu}}({\delta}(\mathbf{Z})\notin \mathcal{A})>0$ for some $\bolds{\nu}\in\Omega$ is inadmissible and dominated by the projection of ${\delta}(\mathbf{Z})$ on $\mathcal{A}$ if either of the following  conditions holds:
\begin{longlist}[(C3)]
  \item[(C3)]  $l(\mathcal{A})\leq (2/3)\pi$;
  \item[(C4)] distribution of ${\delta}(\mathbf{Z})$ is symmetric about ${h}(\bolds{\nu})$ and    with respect to Lebesgue measure, it  satisfies one of the conditions    \textup{(C1)}  and  \textup{(C2)}  with $\psi={h}(\bolds{\nu})$. 
  \end{longlist}
\end{teo}

For the sake of clarity,   the estimation space $\mathcal{A}$ can be called the  action space only  when its Lebesgue measure is less than or equal to $2\pi/3$.
\br
Note that when estimand  ${h}(\bolds{\nu})$  is forced to lie on an arc of semicircle,  $\mathcal{A}$ is strictly a  subset of $\mathbb{T}$.  If ${h}(\bolds{\nu})$ does not    take value on a semicircle,  $\mathcal{A}=\mathbb{T}$ (\cite{b01}, Theorem~7).
\er

\begin{figure}

\includegraphics{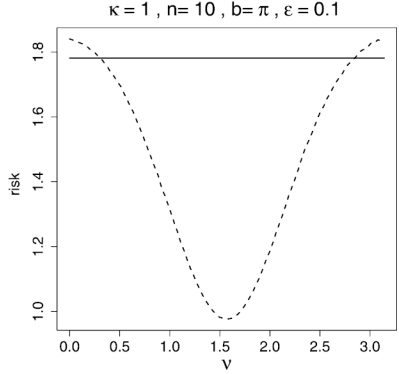}

\caption{Risk plot \textup{(a)} $\bar\theta$ (straight line) and \textup{(b)} projection of $\bar\theta$ on $[0,\pi]$ (dotted line)
under the loss $\mathsf{L}$.}\vspace*{-6pt}\label{fig1}
\end{figure}

Conditions given in  Theorem~\ref{action_space}  for the inadmissibility of  an estimator  are sufficient but not necessary. Suppose that  $\ttt$ has a mixture distribution which is generated from distributions  $\operatorname{CN}(\nu,\kk)$ and $\operatorname{CN}(\nu+\pi,\kk)$ with probabilities $\varepsilon$  and $(1-\varepsilon)$, respectively, where $\nu\in[0,\pi]$, $\varepsilon=0.1$  and $\kk=1$. Based on the  random sample of size $n=10$, risk functions of the sample mean direction $\bar{\ttt}$ (straight line) and the  projection of $\bar{\ttt}$ on $[0,\pi]$ (dotted line) under the loss $\mathsf{L}$  are plotted in Figure~\ref{fig1}. Density of $\operatorname{CN}$  distribution and $\bar{\ttt}$ are defined in the next subsection.  It can be seen that for end points of $\nu\in[0,\pi]$,  $\bar{\ttt}$ is not improved by its projection. This demonstration refutes  the result stated in Theorem~\ref{action_space} for an arbitrary estimator in  case  $l(\mathcal{A})>(2/3)\pi$.

If an estimator  ${\delta}(\mathbf{Z})$ has a  distribution with mixture probability density  $\varepsilon f({\delta}(\mathbf{z})\mid {h}(\bolds{\nu}))+(1-\varepsilon)f({\delta}(\mathbf{z})\mid {h}(\bolds{\nu})+\pi)$, its mean direction is given by
\begin{eqnarray*}
\mathsf{CE}\bigl({\delta}(\mathbf{Z})\bigr)= \cases{ h(\bolds{\nu})+\pi, &\quad
if $\varepsilon<1/2$;
\cr
\mbox{undefined}, &\quad if $\varepsilon=1/2$;
\cr
h(
\bolds{\nu}), &\quad if $\varepsilon>1/2$.}
\end{eqnarray*}
Therefore, when  $\varepsilon<1/2$, ${\delta}(\mathbf{Z})$   can be treated as an  estimator for $ h(\bolds{\nu})+\pi$.

Although Theorem~\ref{action_space} is based on a condition \textup{(C4)}  satisfied by the distribution of an estimator when  $l(\mathcal{A})>(2/3)\pi$, examining the  distribution of the estimator can be  a complex exercise. We try to simplify  these conditions  for specific cases.

Consider the problem of estimating the location parameter $\nu\in\mathbb{T}$  of a circular random variable $\ttt$ under the loss $\mathsf{L}$ which is invariant under a rotation group
\begin{eqnarray}
\label{rotation_G} \mathcal{G}_1= \bigl\{g_\alpha:g_\alpha(
\ttt)=(\ttt+\alpha)\bmod(2\pi) \bigr\}.
\end{eqnarray}
Under  $\mathsf{L}$, an  $\mathcal{G}_1$-equivariant  estimator for $\nu$ based on a random sample $\ttt_1,\ldots,\ttt_n$ satisfies
\begin{eqnarray*}
\delta(\ttt_1,\ldots,\ttt_n)=\ttt_1+\xi(
\ttt_2-\ttt_1,\ldots,\ttt_n-
\ttt_1),
\end{eqnarray*}
where $\xi$ is an arbitrary statistic  whose distribution is free from $\nu$. This  indicates  that   distribution of an $\mathcal{G}_1$-equivariant  estimator  is of the same nature as $\ttt$. Using this fact, we deduce the following result from Theorem~\ref{action_space}.

\begin{cor}\label{mmm1}
If $\ttt$ is a  continuous  circular random variable whose distribution is symmetric about $\nu\in\Omega_1$ and satisfies one of  the conditions    \textup{(C1)}  and  \textup{(C2)}  with $\psi=\nu$ such that $\mathcal{A}=\mathsf{cc}(\Omega_1)$ is $K_1$-type, that is, $l(\mathcal{A})\leq\pi$, any $\mathcal{G}_1$-equivariant  estimator $\delta(\ttt)$ lying outside $\mathcal{A}$ with the positive probability   is inadmissible and dominated by its projection  on $\mathcal{A}$  under loss $\mathsf{L}$.
\end{cor}

A similar result can be extended to the torus $\mathbb{T}^k=\mathbb{T}\times\cdots\times\mathbb{T}$. A distribution on  $\mathbb{T}^k$ can be specified as that of $k$ circular random variables, that is, $k$-tuple  vector  $\mathbf{Z}=(\ttt_1,\ldots,\ttt_k)$ taking values on $\mathfrak{Z}\subset\mathbb{T}^k$.   Suppose that all $k$ components are independently distributed and  each component has a common location parameter $\nu$.  This estimation problem  is invariant under a   group $\mathcal{G}_2$ given by
\begin{eqnarray}
\label{rotation_Grr} \mathcal{G}_2= \bigl\{\mathbf{g}_\alpha=(g_{1\alpha},
\ldots,g_{k\alpha}):g_{i\alpha}(\ttt_i)=(
\ttt_i+\alpha)\bmod(2\pi) \bigr\}.
\end{eqnarray}
 Problem of estimating $\nu$  can also be thought  as multisample problem of estimating common  $\nu$. Therefore, we can also draw random samples of different sizes from different components of $\mathbf{Z}$  as components are independently distributed.   As in  Corollary~\ref{mmm1}, we deduce the following result from Theorem~\ref{action_space}.
\begin{cor}\label{mmm2}
Let all components of  a   random variable $\mathbf{Z}=(\ttt_1,\ldots,\ttt_k)$ taking value on $\mathbb{T}^k$ be independently distributed.
 If  each component has a common location parameter $\nu\in\Omega_1$ and satisfies with respect to Lebesgue measure  one of  the conditions    \textup{(C1)}  and  \textup{(C2)}  with $\psi=\nu$, any $\mathcal{G}_2$-equivariant  estimator $\delta(\mathbf{z})$ lying outside $\mathcal{A}=\mathsf{cc}(\Omega_1)$ with the positive probability   is inadmissible and dominated by its projection  on $\mathcal{A}$  under loss $\mathsf{L}$ when $l(\mathcal{A})\leq\pi$.
\end{cor}

Corollaries~\ref{mmm1} and~\ref{mmm2} enable us to improve various estimators  available in the literature   for the circular location $\nu$ of several  directional  distributions. Apart from the  maximum likelihood estimator (MLE) $\delta_{\mathrm{ml}}$, the following estimators  for $\nu$   have been proposed   on the basis of a random sample $\ttt_1,\ldots,\ttt_n$.
\begin{longlist}[(E2)]
\item[(E1)] (Watson \cite{wbok}, page 135) Sample mean direction $\bar{\ttt}$ minimizes $\sum_{i=1}^n d(\ttt_i,\alpha)$ over  $\alpha\in\mathbb{T}$
  and is obtained as $\bar{\ttt}=\mathsf{atan} (\sum_{i=1}^n\sin\ttt_i/\sum_{i=1}^n\cos\ttt_i )$.
  \item[(E2)] (Mardia and Jupp \cite{mbok}, page~167) Circular median  $\delta_{\mathrm{cm}}$ minimizes $\sum_{i=1}^n d_1(\ttt_i,\alpha)$ over  $\alpha\in\mathbb{T}$, where $d_1(\aaa,\bb)=\pi-
  \llvert  \pi-\llvert  \aaa-\bb\rrvert  \rrvert  $ for $\aaa,\bb\in\mathbb{T}$.
\item[(E3)] (He and Simpson \cite{hes92}) $L_1$-estimator $\delta_{l1}$  minimizes $\sum_{i=1}^n \{d(\ttt_i,\alpha)\}^{1/2}$ over  $\alpha\in\mathbb{T}$.
 \item[(E4)] (Ducharme and Milasevic \cite{dms87a})  Normalized spatial median $\delta_{\mathrm{nsm}}=\mathsf(\alpha_2^*/\alpha_1^*)$, where $(\alpha_1^*,\alpha_2^*)$ is the solution of
\begin{eqnarray*}
\min_{(\alpha_1,\alpha_2)\in\mathbb{R}^2} \sum_{i=1}^n
\bigl\{1+\alpha_1^2+\alpha_2^2-2(
\alpha_1\cos\ttt_i+\alpha_2\sin
\ttt_i)\bigr\}^{1/2}.
\end{eqnarray*}
\item[(E5)] (Neeman and Chang \cite{n01}, Tsai \cite{t09}) Circular Wilcoxon estimator $\delta_{\mathrm{cw}}$  minimizes $\sum_{i=1}^n R_id_1(\ttt_i,\alpha)$ over  $\alpha\in\mathbb{T}$, where $R_i$  $(i=1,\ldots,n)$ is the rank of $\sin(\ttt_i-\alpha)$ amongst $\sin(\ttt_1-\alpha),\ldots,\sin(\ttt_n-\alpha)$.
\end{longlist}
Except $\delta_{\mathrm{cm}}$, all other estimators are proposed in their Cartesian forms. All the above mentioned estimators are either $M$-estimator or  restricted $M$-estimator or $R$-estimator.
\br\label{remz1}
  Due to lack of  preference for the zero direction and orientation in the definition of~$\mathbb{T}$, rotation-equivariant  estimators for circular parameters  are  preferred.  All the above-mentioned estimators enjoy the property of $\mathcal{G}_1$-equivariance  when support is $\mathbb{T}$. When support is  $\mathbb{T}^k$, these estimators based on sample values $\ttt_{11},\ldots,\ttt_{1n_1},\ldots,\ttt_{k1},\ldots,\ttt_{kn_k}$ are also \mbox{$\mathcal{G}_2$-}equivariant. This equivariance property is used in the following section to derive improved estimators.
\er
\subsection{Applications of  Theorem \texorpdfstring{\protect\ref{action_space}}{3.1}}\label{apple1}

Theorem~\ref{action_space} is applicable to a wide variety of estimation problems for directional distributions. In this section, we consider various examples where Theorem~\ref{action_space} leads to improvement over traditional estimators.

\bex[(Unimodal  distributions on circle $\mathbb{T}$)]\label{ex1}
  A circular normal distribution  $\operatorname{CN}(\nu, \kk)$  is defined  by the following density:
\begin{eqnarray*}
f_2(\theta; \nu, \kk) =\frac{1}{2\pi I_{0}(\kk)} e^{\kk \cos (\theta-\nu)},\qquad
\ttt,\nu \in \mathbb{T}, \kk>0,
\end{eqnarray*}
where $I_v$ is the modified Bessel function of the first kind and  order $v$. It is known a priori that $\nu\in\Omega_1$ such that $\Omega_1$ is an arc of the circle.  Without loss of generality, we can assume that $\Omega_1=[0,b]$, where $b\in\mathbb{T}$.  The estimation space for $\nu$ is
\begin{eqnarray*}
\mathcal{A}=\cases{ [0,b], &\quad if $b\leq\pi$;
\cr
\mathbb{T}, &\quad if $b>
\pi$.}
\end{eqnarray*}
The unrestricted MLE of $\nu$ is $\delta_{\mathrm{ml}}=\bar{\ttt}$.  The circular normal distribution is  the only rotationally symmetric distribution for which MLE of the mean direction $\nu$ is the sample mean direction $\bar{\ttt}$. Maximization of likelihood function over  $\Omega_1$ yields the restricted MLE $\delta_{\mathrm{rml}}$  as
\begin{eqnarray}
\label{reses} \delta_{\mathrm{rml}}=\cases{ \bar{\ttt},&\quad if $\bar{\ttt}
\in [0,b]$;
\cr
b, &\quad if $\bar{\ttt}\in (b,\pi+b/{2} )$;
\cr
0, &\quad if $
\bar{ \ttt}\in [\pi+b/{2},2\pi )$;}\quad\mbox{or}\quad \cases{ \bar{\ttt},&
\quad if $ \bar{\ttt} \in [0,b]$;
\cr
b, &\quad if $\bar{\ttt}\in (b,\pi+b/{2}]$;
\cr
0, &\quad if $\bar{\ttt}\in (\pi+b/{2},2\pi )$. }
\end{eqnarray}
At $\bar{\ttt}=\pi+b/{2}$, $\delta_{\mathrm{rml}}$  can take two values.  Since ${\Pr} (\bar{\ttt}=\pi+b/{2} )=0$, both estimators are equivalent.  Note that when $b\leq\pi$, $\delta_{\mathrm{rml}}$ is also the projection of  $\bar{\ttt}$ on $\mathcal{A}$. Corollary~\ref{mmm1} yields that  $\delta_{\mathrm{rml}}$ improves $\bar{\ttt}$ under the loss $\mathsf{L}$ when $b\leq\pi$. When $b>\pi$,  the projection of  $\bar{\ttt}$ on $\mathcal{A}$ is the same as $\bar{\ttt}$.  In a similar way, improvements over all other estimators, $\delta_{\mathrm{cm}}$, $\delta_{l1}$, $\delta_{\mathrm{nsm}}$ and $\delta_{\mathrm{cw}}$ (as defined in \textup{(E2)}, \textup{(E3)}, \textup{(E4)}, \textup{(E5)}) can be obtained from Corollary~\ref{mmm1} and Remark~\ref{remz1}  for   $\nu$ when $\nu$ is restricted to $[0,b]$ and   $b\leq\pi$.

Other well-known symmetric  unimodal distributions  are  wrapped Cauchy $\mathrm{WC}(\nu,\rho)$, wrapped  normal $\operatorname{WN}(\nu,\rho)$  and cardioid $C(\nu,\rho)$  with the following probability  densities   in terms of $t=\cos(\ttt-\nu)$,
\begin{eqnarray*}
f_3(t)&=&\frac{(1-\rho^2)}{(1+\rho^2-2\rho t)},\qquad \rho\in(0,1);
\\
f_4(t)&=&\frac{1}{2\pi}+\frac{1}{\pi}\sum
_{i=1}^\infty \rho^{i^2}T_i(t),
\qquad \rho\in(0,1);
\\
f_5(t)&=&(2\pi)^{-1}+ \pi^{-1}\rho t,\qquad
\llvert \rho\rrvert <1/2,
\end{eqnarray*}
respectively,  where $T_v$ is a Chebyshev polynomial of first kind of order $v$. Jones and Pewsey
\cite{jp05} proposed a family of symmetric unimodal distributions on $\mathbb{T}$ whose densities are provided in terms of $t$ as
\begin{eqnarray*}
f_6(t)=\frac{\{\cosh(\kk\psi)\}^{1/\psi}}{2\pi P_{1/\psi}(\cosh(\kk\psi))} \bigl\{1+\tanh(\kk\psi) t \bigr
\}^{1/\psi},\qquad \kk>0, \psi\in\mathbb{R},
\end{eqnarray*}
 where  $P_v$ is the associated  Legendre   function of the first kind  of degree $v$ and order $0$. Here, we exclude the case of $\kk=0$ since   it yields the uniform distribution on $\mathbb{T}$. Circular normal,   wrapped Cauchy   and cardioid distributions are contained in this family corresponding to  $\psi=0,-1$ and~$1$, respectively.

Another general family of symmetric  unimodal distributions on $\mathbb{T}$ contains wrapped $\alpha$-stable distributions with densities  of the following form  in terms of $t$
(\cite{mbok}, page 52)
\begin{eqnarray*}
f_7(t)&=&\frac{1}{2\pi}+\frac{1}{\pi}\sum
_{i=1}^\infty \rho^{i^\alpha}T_i(t),
\qquad \rho\in(0,1), \alpha=(0,1)\cup(1,2].
\end{eqnarray*}
For $\alpha=2$, it  yields wrapped normal distribution.

Since the distributions corresponding to densities $f_6$ and $f_7$ satisfy condition \textup{(C1)} ,  Corollary~\ref{mmm1} and
 Remark~\ref{remz1} yield improvements over all estimators $\bar{\ttt}$, $\delta_{\mathrm{ml}}$,  $\delta_{\mathrm{cm}}$, $\delta_{l1}$, $\delta_{\mathrm{nsm}}$ and $\delta_{\mathrm{cw}}$ for the  mean direction $\nu$  when $\nu$ is restricted to $[0,b]$ such that $b\leq\pi$.
\eex

\bex[(Mixture distributions on  $\mathbb{T}$)] \label{ex2}
      Let random variable $\ttt$ be generated from $\operatorname{CN}(\nu,\kk)$ and $\operatorname{CN}(\nu+\pi,\kk)$ with probabilities  $\varepsilon$   and $(1-\varepsilon)$, respectively.   For this distribution, $\zeta(t)=e^{2\kk t}$ which is increasing in $t$. Thus, Corollary~\ref{mmm1} and Remark~\ref{remz1} yield that  all the estimators,  $\bar{\ttt}$, $\delta_{\mathrm{ml}}$,  $\delta_{\mathrm{cm}}$, $\delta_{l1}$, $\delta_{\mathrm{nsm}}$ and $\delta_{\mathrm{cw}}$, for $\nu\in\Omega_1$ are improved by their projections on  $\mathcal{A}$ if $l(\mathcal{A})\leq\pi$ and $\varepsilon\geq1/2$. For the values of $\varepsilon<1/2$, improvements are possible  when  $l(\mathcal{A})\leq(2/3)\pi$.

 Similar improvements are possible    when we  mix either two wrapped normal or two  distributions with density $f_6$  for same value of $\psi\leq1$ with different mean directions $\nu$ and $\nu+\pi$. To apply  Corollary~\ref{mmm1}, we must  show that the corresponding  function $\zeta(t)$ is an increasing function in $t$.  For  a general density $f_6(t)$, derivative of $\zeta(t)$  with respect to $t$ is obtained as
\begin{eqnarray*}
\zeta'(t)=\frac{2(1-\psi)}{\psi}\frac{\tanh(\kk\psi)}{\{1+\tanh(\kk\psi)  t\}^2} \biggl\{
\frac{1+\tanh(\kk\psi)  t}{1-\tanh(\kk\psi)  t} \biggr\}^{1/\psi},\qquad \psi\neq0,
\end{eqnarray*}
which is always nonnegative unless $\psi>1$.

The wrapped normal distribution can be represented by theta function $\vartheta_3$. Using the representation of  $\vartheta_3$  in terms of infinite products
(\cite{gr65}, page~921, equation~8.181.2), we have
\begin{eqnarray*}
f_4(t)=\frac{1}{2\pi}\prod_{i=1}^{\infty}
\bigl(1-\rho^{2i}\bigr) \bigl(1+2t\rho^{2i-1}+
\rho^{2(i-1)} \bigr).
\end{eqnarray*}
Second derivative of $f_4(t)$ with respect to $t$ is $f_4(t)\sum_{i\neq j}\xi_i(t)\xi_j(t)$, where $\xi(i)={2\rho^{2i-1}}/\{1+2t\rho^{2i-1}+\rho^{2(i-1)}\}^{-1}$. Convexity of $f_4(t)$ follows from the positiveness of $\xi_i(t)$ and increasing nature of $\zeta(t)$ follows from the convexity of $f_4(t)$ using Remark~\ref{zee1}.
\eex

\bex[(Distributions on $\mathbb{T}$ with $k$-fold rotational symmetry)] \label{ex3}
  This distribution is constructed by putting $k$  copies of the original distribution end-to-end
  (\cite{mbok}, page~53). If we are given a   distribution of $\ttt$ which is  unimodal and symmetric  about $\nu$, constructed distribution has the density  $f(\cos(k(\ttt-\nu_0)))$,   $\ttt\in\mathbb{T}, \nu_0\in[0,2\pi/k)$. Note that new  distribution is  $k$-modal, for example, $k$-modal circular normal distribution
(\cite{js01}, page~209). For $k\geq3$, condition \textup{(C3)}  of Theorem~\ref{action_space} is satisfied and so estimators of $\nu_0$ lying outside $\mathcal{A}$ with a positive probability can be improved by  their projections on $\mathcal{A}$. Note that the results hold when the parameter space is full, that is, $[0,2\pi/k)$ or restricted, that is, a subset of  $[0,2\pi/k)$. For $k=2$, the result holds only for restricted parameter space if   $\nu_0\in\Omega_1\subset[0,b]$ such that $b\leq(2/3)\pi$.
  \eex

\bex[(Distributions on torus  $\mathbb{T}^k$)] \label{ex4}
   Suppose that all $k$ components of $\mathbf{Z}=(\ttt_1,\ldots,\ttt_k)\in\mathfrak{Z}\subset\mathbb{T}^k$ are independently distributed and  $i$th component $\ttt_i$ follows $\operatorname{CN}(\nu,\kk_i)$ with known $\kk_i$. Consider the estimation of common $\nu\in\Omega_1$ under the loss $\mathsf{L}$. Let   $(\ttt_{i1},\ldots,\ttt_{in_i})$ be a random sample from the $i$th population $(i=1,\ldots,n)$.  Suppose $\bar{\ttt}_i$ denotes the sample mean direction and $R_i=\{(\sum_{j=1}^{n_i}\sin\ttt_{ij})^2+(\sum_{j=1}^{n_i}\cos\ttt_{ij})^2\}^{1/2}$ denotes the sample resultant length for the sample   of $i$th component. The  MLE of   $\nu$ is
\begin{eqnarray}
\label{to_e} \tilde{\ttt}=\mathsf{atan} \biggl(\frac{\sum_{i=1}^{k}\kk_iR_i\sin\bar{\ttt}_i}{\sum_{i=1}^{k}\kk_iR_i\cos\bar{\ttt}_i} \biggr).
\end{eqnarray}
 The conditional distribution of $\tilde{\ttt}$ is again circular normal $\operatorname{CN}(\nu,R^*)$,  where
\[
R^*=\Biggl\{\Biggl(\sum_{i=1}^{k}\kk_iR_i\sin\bar{\ttt}_i\Biggr)^2+\Biggl(\sum_{i=1}^{k}\kk_iR_i\cos\bar{\ttt}_i\Biggr)^2  \Biggr\}^{1/2}\qquad
\mbox{Holmquist \cite{h91}.}
\]
Since the distribution of $R^*$ is dependent only on $\kk_i$, distribution  of   $\tilde{\ttt}$ is unimodal and symmetric about its mode $\nu$. Corollary~\ref{mmm2} yields that   $\tilde{\ttt}$ is dominated by its projection on $\mathcal{A}$ under the loss $\mathsf{L}$ if $l(\mathcal{A})\leq\pi$. Note that improved estimator of $\tilde{\ttt}$  is also the restricted MLE for common $\nu$ when $l(\mathcal{A})\leq\pi$.

 Using Corollary~\ref{mmm2}  and Remark~\ref{remz1}, we can obtain improvements over other estimators  $\bar{\ttt}$, $\delta_{\mathrm{cm}}$, $\delta_{l1}$, $\delta_{\mathrm{nsm}}$ and $\delta_{\mathrm{cw}}$ based on random sample $\ttt_{11},\ldots,\ttt_{1n_1},\ldots,\ttt_{k1},\ldots,\ttt_{kn_k}$.
\br
This  model  can be further extended to the cases where distributions of independent components are not necessarily the same and satisfy conditions \textup{(C1)}  and  \textup{(C2)}  with respect to Lebesgue measure, namely, distributions considered in   Examples~\ref{ex1} and~\ref{ex2}.
 \er
\eex

\bex[(Distribution on unit sphere $\mathbb{S}_{2}$)] Point $\mathbf{x}=(x_1,x_2,x_3)^T\in\mathbb{S}_{2}$ can be specified by its geographical coordinates: colatitude  $\ttt=\cos^{-1}(x_3)\in[0,\pi]$ and longitude  $\phi=\mathsf{atan}(x_2/x_1)\in[0,2\pi)$. Hence,  sphere $\mathbb{S}_{2}$ is isomorphic to $\mathbb{T}_2=[0,\pi]\times[0,2\pi)$.   The density of Fisher distribution on the support $\mathbb{T}_2$ is given by
\begin{eqnarray}
f_8(\theta, \phi; \nu_{1},\nu_{2}, \kk) =
\frac{\kk \sin\ttt}{4\pi \sinh (\kk)} \exp\bigl\{\kk \bigl(\sin\ttt \sin\nu_{1}\cos (\phi-
\nu_{2})+\cos \ttt\cos\nu_{1}\bigr)\bigr\},
\nonumber
\\
\eqntext{(\nu_{1},\nu_{2})\in\mathbb{T}_{2},\kk>0,}
\end{eqnarray}
where $(\nu_{1},\nu_{2})$ is the mean direction. On the basis of a random sample of size $n$, the  MLE of $\nu_2$, when $\kk$ is known/unknown and $\nu_1$ is unknown,  is given by
\begin{eqnarray*}
\delta_{\mathrm{ml}}=\mathsf{atan} \biggl(\frac{\sum_{i=1}^n \sin\ttt_i\sin\phi_i}{\sum_{i=1}^n \sin\ttt_i\cos\phi_i} \biggr).
\end{eqnarray*}
It is known that mean direction is restricted to a continuous arc of hemisphere.  Without loss of generality, we can assume that $\nu_1\in[0,\pi]$ and $\nu_{2}\in[0,b]$ with $b\leq\pi$. We can improve the MLE of $\nu_2$ by its projection on $[0,b]$  if $b\leq(2/3)\pi$ using Theorem~\ref{action_space}.
\eex

\bex[(Distribution on cylinder $\mathbb{R}\times\mathbb{T}$)]
Mardia and Sutton
\cite{m78} proposed a distribution  on the cylinder $\mathbb{R}\times\mathbb{T}$. Let $(X,\ttt)$ have the   support $\mathbb{R}\times\mathbb{T}$ where  the marginal distribution of $\ttt$ is $\operatorname{CN}(\nu_0,\kk)$, $\nu_0\in\mathbb{T}, \kk>0$ and conditional  distribution of $X\mid \ttt$ is a normal  with mean
     \begin{eqnarray*}
\mu_c=\mu+\sigma\rho\sqrt{\kk}\bigl\{\cos(\ttt-\nu)- \cos(
\nu_0-\nu)\bigr\},\qquad \mu\in\mathbb{R}, \nu\in\mathbb{T},0\leq\rho
\leq1, \sigma>0
\end{eqnarray*}
 and variance $\sigma^2(1-\rho^2)$. Based on a random sample of size $n$,  the MLE of $\nu$ is
 \begin{eqnarray*}
\delta_{\mathrm{ml}}=\mathsf{atan} \biggl(\frac{s_2}{s_3} \frac{r_{23}r_{12}-r_{13}}{r_{23}r_{13}-r_{12}}
\biggr),
\end{eqnarray*}
 where for $i=1,2,3$ and $j=1,\ldots,n$, $x_{1j}=x_j$, $ x_{2j}=\cos\ttt_j$, $ x_{3j}=\sin\ttt_j$, $\bar{x}_i=\sum_{j}x_{ij}/n$, $ s_i^2=\sum_{j}(x_{ij}-\bar{x}_i)^2$;  and for $i\neq k$ $(i,k=1,2,3)$,
  \begin{eqnarray*}
r_{ik}&=&\frac{1}{s_is_k} \sum_{j}
(x_{ij}-\bar{x}_i) (x_{kj}-\bar{x}_k).
\end{eqnarray*}
If $\nu$ is restricted to  $\Omega_1 $ such that  $l(\mathcal{A})\leq (2/3)\pi$,  $\delta_{\mathrm{ml}}$ is dominated by its
projection on $\mathcal{A}$ under the loss $\mathsf{L}$ using Theorem~\ref{action_space}. Note that simulations indicate that this result is also valid for  $(2/3)\pi<l(\mathcal{A})\leq\pi$.
\eex

\section{An inadmissibility result for  general equivariant rules}\label{SE4}

Moors \cite{moors81,moors85}, Kumar and Sharma
\cite{ks92} gave a general method for obtaining improved equivariant estimators of parameters in Euclidean spaces. In this section, we extend these results for estimating circular parameters.

Let the   problem of estimating $h(\bolds{\nu})\in  \Omega_1$ under the loss $\mathsf{L}$ be invariant under a compact group  $\mathcal{G}$ of measurable transformations $\mathbf{g}:\mathfrak{Z}\rightarrow\mathfrak{Z}$. There exists  a  finite and  left (right)  invariant Haar measure $\lambda$ on $\mathcal{G}$
(\cite{e89}, Theorem 1.5). Let $\bar{ \mathcal{G}}$ and $\tilde{\mathcal{G}}$ be the groups induced  by $\mathcal{G}$ on parameter space $\Omega$ and estimation space $\mathcal{A}=\mathsf{cc}(\Omega_1)$.

\begin{lem}\label{lem}
For the  $\mathcal{G}$-invariant estimation problem defined above, we have:
\begin{longlist}[(iii)]
\item[(i)] $f(\mathbf{z}\mid\bar{\mathbf{g}}(\bolds{\nu}))=f(\mathbf{g}^{-1}(\mathbf{z})\mid \bolds{\nu})$ a.e. with respect to measure $\eta$;
  \item[(ii)] ${h}\bar{\mathbf{g}}(\bolds{\nu})=\tilde{{g}}{h}(\bolds{\nu})\bmod(2\pi)$  for all  $\bolds{\nu}\in\Omega$;
  \item[(iii)] Let $A\in\mathfrak{C}_5$ (closed convex). If  ${\phi}_{0}$ is the projection of ${\phi}\in\mathbb{T}$ on $A$, the projection of $\tilde{g}(\phi)$ on $\tilde{{g}}(A)$ is $\tilde{{g}}(\phi_0)$.
 \end{longlist}
\end{lem}

For each $\mathbf{z}\in\mathfrak{Z}$ and $\bolds{\nu}\in\Omega$,  define a probability measure on $\mathcal{G}$ as
\begin{eqnarray*}
\tau \bigl(\mathbf{z}\mid \bar{\mathbf{g}}(\bolds{\nu}) \bigr)=
\frac{f (\mathbf{z}\mid \bar{\mathbf{g}}(\bolds{\nu}) )}{\int_{\mathcal{G}} f (\mathbf{z}\mid \bar{\mathbf{g}}^*(\bolds{\nu}) ) \, \mathrm{d}\lambda(\mathbf{g}^*)}.
\end{eqnarray*}
With the help of these measures, for a fixed $\mathbf{z}\in\mathfrak{Z}$, define a function ${h}_{\mathbf{z}}:\Omega\rightarrow \mathcal{A}$ as
\begin{eqnarray}
\label{h} {h}_{\mathbf{z}}(\bolds{\nu}) =\cases{\displaystyle \mathsf{atan} \biggl(
\frac{ \int_{\mathcal{G}} \sin  (\tilde{g}{h}(\bolds{\nu}) ) \tau (\mathbf{z}\mid \bar{\mathbf{g}}(\bolds{\nu}) )  \,\mathrm{d}\lambda(\mathbf{g})}{\int_{\mathcal{G}} \cos  (\tilde{g}{h}(\bolds{\nu})  ) \tau (\mathbf{z}\mid \bar{\mathbf{g}}(\bolds{\nu}) )  \,\mathrm{d}\lambda(\mathbf{g})} \biggr),&\quad if $\displaystyle\int_{\mathcal{G}} f \bigl(
\mathbf{z}\mid \bar{\mathbf{g}}(\bolds{\nu}) \bigr) \,\mathrm{d} \lambda(
\mathbf{g})>0$;
\vspace*{3pt}\cr
{h}(\bolds{\nu}),&\quad if $\displaystyle\int_{\mathcal{G}} f
\bigl(\mathbf{z}\mid \bar{\mathbf{g}}( \bolds{\nu}) \bigr) \,\mathrm{d}\lambda(
\mathbf{g})=0$.}
\end{eqnarray}
The new estimation space  $\mathcal{A}_{\mathbf{z}}$ is defined as a convex closure of ${h}_{\mathbf{z}}(\Omega)$, that is, $\mathcal{A}_{\mathbf{z}}=\mathsf{cc} ( {h}_{\mathbf{z}}(\Omega) )$.
\br\label{subset}
In  case of $\int_{\mathcal{G}} f (\mathbf{z}\mid \bar{\mathbf{g}}(\bolds{\nu}) ) \,\mathrm{d}\lambda(\mathbf{g})=0$, $\mathcal{A}_{\mathbf{z}}=\mathcal{A}$. When $\int_{\mathcal{G}} f (\mathbf{z}\mid \bar{\mathbf{g}}(\bolds{\nu}) ) \,\mathrm{d}\lambda(\mathbf{g})>0$,  ${h}_{\mathbf{z}}(\bolds{\nu}) $ can be written as
\begin{eqnarray}
\label{res_h} {h}_{\mathbf{z}}(\bolds{\nu}) =\mathsf{atan} \biggl(
\frac{\mathsf{E} \sin  (\tilde{g}{h}(\bolds{\nu}) ) }{\mathsf{E}  \cos  (\tilde{g}{h}(\bolds{\nu}) )} \biggr),
\end{eqnarray}
where the expectation is taken over $\mathbf{g}$ with respect to a probability measure $\tau (\mathbf{z}\mid \bar{\mathbf{g}}(\bolds{\nu}) )  \,\mathrm{d}\lambda(\mathbf{g})$. If $\Pr(\tilde{g}{h}(\bolds{\nu})\in\mathcal{A})=1$,    ${h}_{\mathbf{z}}(\bolds{\nu})\in \mathcal{A}$ from the convexity of $\mathcal{A}$ (see  Proposition~\ref{jen}). Since $\mathcal{A}_{\mathbf{z}}$ is the smallest  convex set containing ${h}_{\mathbf{z}}(\bolds{\nu})$, we conclude that  $\mathcal{A}_{\mathbf{z} }\subset\mathcal{A}$.
\er

We assume that every transformation  $\tilde{{g}}\in\tilde{\mathcal{G}}$ satisfies  circular property, that is, for a fixed $\alpha$ and for all $\phi\in\mathcal{A}$, $\tilde{{g}}(\phi)$ is either $\aaa+\phi$   or $\aaa-\phi$ such that images are in  $\mathcal{A}$ itself.

\begin{lem}\label{ac}
If  every induced transformation   on  $\mathcal{A}$ satisfies the circular property, the reduced estimation space  $\mathcal{A}_{\mathbf{z} }$ satisfies $\mathcal{A}_{\mathbf{g} (\mathbf{z})}=\tilde{{g} } (\mathcal{A}_{\mathbf{z} })$ for all $\mathbf{z}\in\mathfrak{Z}$ and $\mathbf{g} \in \mathcal{G}$.
\end{lem}

Lemmas~\ref{lem} and~\ref{ac} are utilized to prove the  main result  of this section.

\begin{teo}\label{main}
Consider an $\mathcal{G}$-invariant estimation problem under the loss function $\mathsf{L}$ with a compact group $\mathcal{G}$ such that  elements of  induced group $\tilde{\mathcal{G}}$ on the  estimation space $\mathcal{A}$ satisfy  the circular property. Any $\mathcal{G}$-equivariant estimator ${\delta}$ satisfying ${\Pr}_{\bolds{\nu}}({\delta}(\mathbf{Z})\notin\mathcal{A}_{\mathbf{Z}})>0$  for some $\bolds{\nu}\in \Omega$ is dominated  by its projection  on $\mathcal{A}_{\mathbf{Z}}$ provided that $l(\mathcal{A}_{\mathbf{Z}})\leq (2/3)\pi$.
\end{teo}

This result is applicable to both restricted and unrestricted estimation problems as illustrated in  the following subsection.

\subsection{Applications of Theorem \texorpdfstring{\protect\ref{main}}{4.1}}\label{app}
We consider  estimation of the  location parameter $\nu\in\Omega_1$ of a   circular random variable $\ttt$  under the loss $\mathsf{L}$. Let us denote by $f_{\nu}(\ttt)$ the density of $\ttt$ that would be a function of $\cos(\ttt-\nu)$.

\subsubsection{Unrestricted estimation problems}\label{urstt}

The estimation problem is invariant under the  rotation group $\mathcal{G}_1$. Clearly, the induced group $\tilde{\mathcal{G}}_1$ on the estimation space $\mathcal{A}=\mathbb{T}$ is itself  $\mathcal{G}_1$  and every transformation in $\mathcal{G}_1$ satisfies  the circular  property. Taking $\mathcal{G}=\mathcal{G}_1$,   we define
\begin{eqnarray*}
h_{\ttt}(\nu)=\mathsf{atan} \biggl(\frac{\int_{0}^{2\pi}\sin(\alpha+\nu) f_{\aaa+\nu}(\ttt) \,{\rm{d}} \alpha}{\int_{0}^{2\pi}\cos(\alpha+\nu) f_{\aaa+\nu}(\ttt) \,{\rm{d}} \alpha} \biggr)=
\mathsf{atan} \biggl(\frac{\int_{0}^{2\pi}\sin\alpha  f_{\aaa}(\ttt) \,{\rm{d}} \alpha}{\int_{0}^{2\pi}\cos\alpha f_{\aaa}(\ttt) \,{\rm{d}} \alpha} \biggr),
\end{eqnarray*}
or equivalently,  $h_{\ttt}(\nu)$ is constant on $\nu\in\mathbb{T}$. Therefore  the following result  follows from Theorem~\ref{main}.

\begin{cor}\label{ccoo}
If $\ttt$ is a circular random variable with the unrestricted location parameter $\nu\in\mathbb{T}$, there is only one admissible $\mathcal{G}_1$-equivariant  estimator under the loss $\mathsf{L}$ which is obtained as
\begin{eqnarray*}
\delta_{\mathrm{ad}}=\mathsf{atan} \biggl(\frac{\int_{0}^{2\pi}\sin\alpha  f_{\aaa}(\ttt) \,{\rm{d}} \alpha}{\int_{0}^{2\pi}\cos\alpha f_{\aaa}(\ttt) \,{\rm{d}} \alpha} \biggr),
\end{eqnarray*}
where  $f_{\nu}(\ttt)$ the density of $\ttt$.
\end{cor}

For $\operatorname{CN}(\nu,\kk)$ distribution,   $\delta_{\mathrm{ad}}=\bar{\ttt}$. All other $\mathcal{G}_1$-equivariant  estimators $\delta_{\mathrm{cm}}$, $\delta_{l1}$, $\delta_{\mathrm{nsm}}$ and $\delta_{\mathrm{cw}}$ for  $\nu\in\mathbb{T}$ under the loss $\mathsf{L}$ are improved by $\bar{\ttt}$ using Corollary~\ref{ccoo}. This result is significant in the sense that so far comparison of $\bar{\ttt}$ with $\delta_{\mathrm{cm}}$, $\delta_{l1}$, $\delta_{\mathrm{nsm}}$ and $\delta_{\mathrm{cw}}$ was done only with respect to asymptotic efficiency and robustness.

\subsubsection{Restricted estimation problems}\label{rstt}
Let the location   $\nu$ be  restricted to any arc of semicircle. Without loss of generality, we can assume that $\nu\in \Omega_1=[0,b]$ with $0<b\leq\pi$. In Example~\ref{ex1}, the same restricted space estimation problem was considered  and we obtained  improvements over estimators $\bar{\ttt}$, $\delta_{\mathrm{ml}}$, $\delta_{\mathrm{cm}}$, $\delta_{l1}$, $\delta_{\mathrm{nsm}}$ and $\delta_{\mathrm{cw}}$ for various distributions.  Denote by  $\bar{\ttt}^*$,  $\delta^*_{\mathrm{ml}}$, $\delta^*_{\mathrm{cm}}$, $\delta^*_{l1}$, $\delta^*_{\mathrm{nsm}}$ and $\delta^*_{\mathrm{cw}}$, the dominating estimators of  $\bar{\ttt}$,  $\delta_{\mathrm{ml}}$, $\delta_{\mathrm{cm}}$, $\delta_{l1}$, $\delta_{\mathrm{nsm}}$ and $\delta_{\mathrm{cw}}$ as their projections on $\Omega_1=[0,b]$.    Here, we can further  improve upon  these improved estimators.

 We consider the transformation group as $\mathcal{G}_3=\{e,g\}$ where $g(\ttt)=(b-\ttt)\bmod(2\pi)$ and $e$ is the identity transformation. The estimation problem remains  invariant under  $\mathcal{G}_3$ and the induced group on the estimation space $\mathcal{A}$ is $\mathcal{G}_3$.  Clearly, elements of   $\mathcal{G}_3$  satisfy the  circular property.  Define the following   function for a fixed  random sample $\ttt_1,\ldots,\ttt_n$  as
\begin{eqnarray*}
h_{({\ttt}_1,\ldots,{\ttt}_n)}(\nu) &=& \bigl\{{b}/{2}+\tan^{-1} \bigl(a(\nu) \bigr)
\bigr\}\bmod(2\pi),
\end{eqnarray*}
where after some algebraic computations, $a(\nu)$ is derived  as
\begin{eqnarray*}
a(\nu)&=&\frac{  \prod_{i=1}^{n}   f_{b-\nu}(\ttt_{i})- \prod_{i=1}^{n}  f_{\nu}(\ttt_{i})}{ \prod_{i=1}^{n}  f_{b-\nu}(\ttt_{i}) +
  \prod_{i=1}^{n}   f_{\nu}(\ttt_{i})}\tan \biggl(\frac{b}{2}-\nu \biggr).
\end{eqnarray*}
Since $f_{\nu}(\ttt)$ is a function of  $\cos(\ttt-\nu)$, $h_{({\ttt}_1,\ldots,{\ttt}_n)}(\nu)$  is symmetric about $\nu=b/2$. It is sufficient to assume that $\nu\in[0,b/2]$ to study the behaviour of the function  $h_{({\ttt}_1,\ldots,{\ttt}_n)}(\nu)$.

We consider     distributions for which $a(\nu)$ is monotonic in $\nu\in[0,b/2]$. As $a(b/2)=0$, monotonic nature of $a(\nu)$  is dependent on the sign of $a(0)$. The new  estimation space    is
\begin{eqnarray*}
\mathcal{A}_{(\ttt_1,\ldots,\ttt_n)}=\mathsf{cc} \bigl(h_{({\ttt}_1,\ldots,{\ttt}_n)}(
\Omega_1) \bigr)= \cases{ \bigl[b^{*}, {b}/{2} \bigr],&
\quad if $a(0)<0$;
\cr
\{{b}/{2}\}, &\quad if $a(0)=0$;
\cr
\bigl[{b}/{2},
b^{*} \bigr], &\quad if $a(0)>0$,}
\end{eqnarray*}
where $b^{*}=h_{({\ttt}_1,\ldots,{\ttt}_n)}(0)$. It may be noted that $l(\mathcal{A}_{(\ttt_1,\ldots,\ttt_n)})\leq b/2$ which is a substantial  reduction in $l(\mathcal{A})=b$.

Since $\mathcal{G}_3\subset\mathcal{G}_1$, all estimators $\bar{\ttt}$, $\delta_{\mathrm{ml}}$,  $\delta_{\mathrm{cm}}$, $\delta_{l1}$, $\delta_{\mathrm{nsm}}$, $\delta_{\mathrm{cw}}$, $\bar{\ttt}^*$, $\delta_{\mathrm{ml}}^*$, $\delta^*_{\mathrm{cm}}$, $\delta^*_{l1}$, $\delta^*_{\mathrm{nsm}}$ and $\delta^*_{\mathrm{cw}}$ are also $\mathcal{G}_3$-invariant. For the distributions which satisfy  the assumption of monotonicity of $a(\nu)$, all these estimators can be  improved by their projections on $\mathcal{A}_{(\ttt_1,\ldots,\ttt_n)}$ using Theorem~\ref{main}.

If $\ttt$ follows a  $\operatorname{CN}(\nu,\kk)$ distribution, the function $a(\nu)$ is given as
\begin{eqnarray*}
a(\nu)=\tanh \biggl(\kk r \sin \biggl(\tilde{\ttt}-\frac{b}{2} \biggr)
\sin \biggl(\frac{b}{2}-\nu \biggr) \biggr)\tan \biggl(\frac{b}{2}-
\nu \biggr).
\end{eqnarray*}
Monotonicity of $a(\nu)$ can be easily observed. The new estimation space is equal to
\begin{eqnarray*}
\mathcal{A}_{(\bar{\ttt}, r)}= \cases{ \bigl[b^{*},{b}/{2} \bigr],&
\quad if $\bar{\ttt} \in J_1$;
\cr
\{{b}/{2}\}, &\quad if $\bar{\ttt}
\in\{{b}/{2}$, $\pi+{b}/{2} \}$;
\cr
\bigl[{b}/{2}$, $b^{*} \bigr], &
\quad if $\bar{\ttt}\in I_1$,}
\end{eqnarray*}
where $b^{*}={b}/{2}+\tan^{-1} [\tan({b}/{2})\tanh \{\kk r \sin (\bar{\ttt}-{b}/{2} )\sin({b}/{2})   \} ]$, $J_1= J ({b}/{2},\pi+{b}/{2} )$ and $I_1= I ({b}/{2},\pi+{b}/{2} )$. Based on $\mathcal{A}_{(\bar{\ttt}, R)}$, an estimator  $\delta$ is dominated by
\begin{eqnarray*}
\delta_I=\cases{ \delta, &\quad if $\bar{\ttt}\in J_1$
and $\delta\in \bigl[b^*,{b}/{2} \bigr]$, or, $\bar{\ttt}\in I_1$ and
$\delta\in \bigl[{b}/{2},b^* \bigr]$;
\cr
b^*, &\quad if $\bar{\ttt}\in
J_1$ and $\delta\in J \bigl(b^*,\gamma \bigr)$, or, $\bar{\ttt}\in
I_1$ and $\delta\in \bigl(b^*,\gamma \bigr)$;
\cr
b/2, &\quad if $
\bar{\ttt}\in J_1$ and $\delta\in ({b}/{2},\gamma ]$, or, $\bar{\ttt}
\in I_1$ and $\delta\in J ({b}/{2},\gamma ]$, or, $\bar{\ttt}\in
\{b/2, \pi+{b}/{2}\}$,}
\end{eqnarray*}
where $\gamma=\pi+b^*/2+b/4$. In Example~\ref{ex1}, the estimator $\delta$ for $\nu$ is  dominated by the projection of  $\delta$  on $\mathcal{A}$. Let us denote by   $\delta^*$  this improved estimator.  Based on $\mathcal{A}_{(\bar{\ttt},R)}$, improved estimators  $\delta_I$ and $\delta_I^*$ of $\delta$ and $\delta^*$, respectively,  are equivalent except when either $\delta\in [\gamma,\pi+b/2]$ and  $\bar{\ttt}\in J_1$ or $\delta\in[\pi+b/2,\gamma]$ and $\bar{\ttt}\in I_1$. If $\delta=\bar{\ttt}$,  both $\bar{\ttt}_{I}$ and $\bar{\ttt}_{I}^{*}$ are equivalent and given by
\begin{eqnarray}
\label{iitt} \cases{
 b^{*}, & \quad if $\cases{
     \bar{\ttt}<b^{*}\mbox{ and }\bar{\ttt}\in [0,{b}/{2}),\vspace*{3pt}\cr
     \bar{\ttt}>b^{*}\mbox{ and }\bar{\ttt}\in ({b}/{2},\pi+{b}/{2}),\vspace*{3pt}\cr
     \bar{\ttt}\in (\pi+{b}/{2},2\pi),}$
     \vspace*{3pt}\cr
{b}/{2}, &\quad if $\bar{\ttt}= \pi+{b}/{2}$,\vspace*{3pt}\cr
\bar{\ttt}, &\quad elsewhere.}
\end{eqnarray}
In Figure~\ref{trun11}, we have plotted the risk functions of the MLE $\bar{\ttt}$, restricted estimator $\delta_{\mathrm{rml}}$ defined in (\ref{reses})   and improved estimator  $\bar{\ttt}_{I}$ defined in~(\ref{iitt}) under the loss function $\mathsf{L}$. The risk values have been evaluated using simulations. For this, we have generated 100\,000 samples from a $\operatorname{CN}(\nu, \kk)$ distribution for various values of $(n,\kk,b)$. The following conclusions can be made from the numerical study.
\begin{longlist}[(a)]
  \item[(a)] The risk function of $\bar{\ttt}$ is constant for a fixed value of $(n,\kk)$. The risk functions of $\delta_{\mathrm{rml}}$ and  $\bar{\ttt}_{I}$ are symmetric about $b/2$. For small values of $\kk$ or $b$, these risk functions are strictly decreasing in $\nu\in[0,b/2]$. For higher values of $\kk$ or $b$, behaviour is reverse.

\begin{figure}

\includegraphics{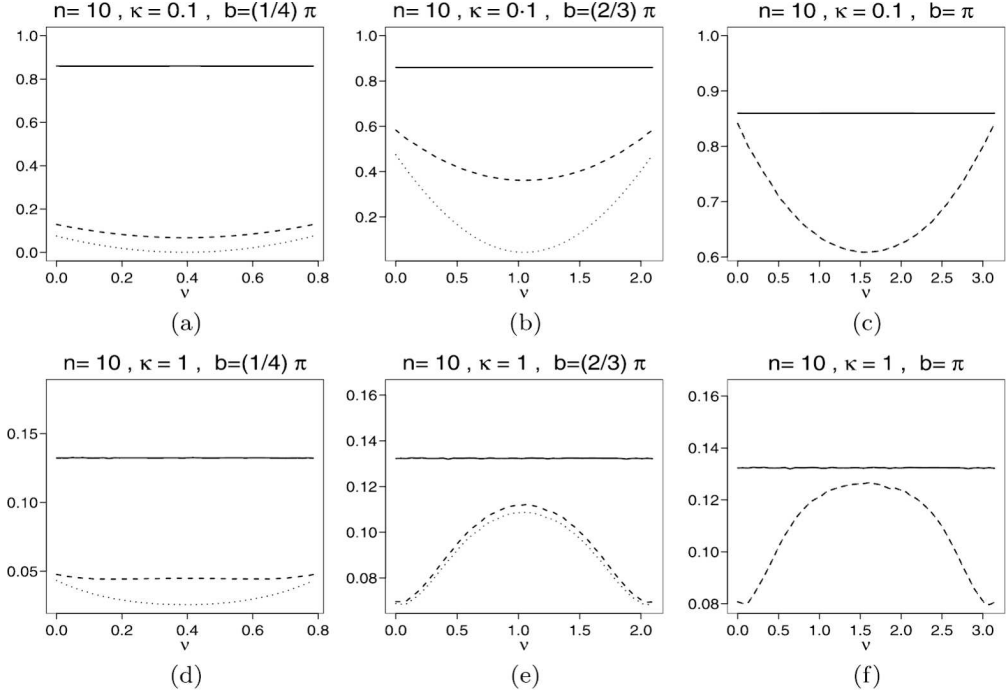}

\caption{Risk plots of (a) $\bar{\ttt}$ (straight line), (b) $\delta_{\mathrm{rml}}$ (dashed line) and (c) $\bar{\ttt}_{I}$ (dotted line) under the loss function $\mathsf{L}$ when $b\in(0, \pi]$.}\label{trun11}
\end{figure}

  \item[(b)]  For all the values of $(n,\kk)$ and $\nu\in[0,b]$, $\delta_{\mathrm{rml}}$ uniformly improves $\bar{\ttt}$   and $\bar{\ttt}_{I}$  uniformly improves  $\delta_{\mathrm{rml}}$ when $b \in (0,\pi]$. Risk values of  $\delta_{\mathrm{rml}}$ and  $\bar{\ttt}_{I}$ are the same when $b=\pi$ and are less than that of $\bar{\ttt}$ for all $\nu\in[0,b]$.

\item[(c)] The amount of relative improvement of    $\delta_{\mathrm{rml}}$   over     $\bar{\ttt}$  is increasing as $\kk$ or $b$ decreases. This is seen to be as high as 95\%. Similarly,  relative improvement of    $\bar{\ttt}_{I}$ over  $\delta_{\mathrm{rml}}$ is seen to be up to~75\%.
\end{longlist}
Similar observations have been made for various other configurations of  $(n,\kk,b)$. Simulations for various directional distributions show significant improvements. We omit details here.

\br
For the support $\mathbb{T}$, both unrestricted and restricted estimation problems discussed in Sections~\ref{urstt} and~\ref{rstt} can be easily extendable to the case of support $\mathbb{T}^k$.
\er

\section{Concluding remarks}

For estimating parameters in Euclidean spaces, with respect to the loss function as an increasing function of distance,  the action space is taken to be the smallest convex set containing the estimand space. If an estimator lies outside it with a positive probability, an improvement is obtained by projecting this estimator on the action space. In Section~\ref{SE3}, we have extended this concept to the estimation of circular parameters of directional distributions. The result of Theorem~\ref{action_space}, is not exactly analogous to the result for Euclidean spaces.  Further, \cite{moors81,ks92} gave a new technique for improving equivariant estimators in Euclidean spaces. In  Section~\ref{SE4}, we have developed a theory to extend this technique to circular parameters. The results have been applied to various estimation problems in directional distributions. The resulting estimators are seen to show significant improvements over the usual estimators.

It would be interesting to further extend these results to parameters lying in spheres of higher dimensions.

\begin{appendix}
\section*{Appendix}
\subsection{Proof of Proposition  \texorpdfstring{\protect\ref{convex_pr}}{2.1}}
The statement trivially follows when  $A=\mathbb{T}$. From Lemma~\ref{rot}, it is sufficient to consider  $A$ to be of $I$-type. First, we assume $A=I[0,\gamma]$ with $\gamma\in(0,\pi]$. Note that for $\gamma=0$, the proof  is  trivial.

Let $\{\phi_1,\ldots,\phi_n\}$ be  linearly ordered sample of  $n(\geq 2)$ observations taking values in $A$. After rotating  sample by an angle $2\pi-\phi_1$, modified ordered sample $\{0,\alpha_2,\ldots,\alpha_{n}\}$ takes values in $I[0,\alpha_{n}]\subset A$ with $\alpha_{i}=\phi_{i}-\phi_1$ $(i=2,\ldots,n)$. We have to prove
\begin{eqnarray*}
\bar{\alpha}_w=\mathsf{atan} \biggl(\frac{\sum_{i=2}^n w_i \sin\alpha_i}{ w_1+ \sum_{i=2}^n w_i \cos\alpha_i} \biggr)\in
I[0,\alpha_n],
\end{eqnarray*}
 where weights are  $w_1,\ldots,w_n$. From (\ref{atan}), the lower bound of $\bar{\alpha}_w$ is zero. If $\alpha_n\neq\pi$,  both $\sum_{i=2}^{n}w_i\sin\alpha_i$ and  $w_1+\sum_{i=2}^{n}w_i\cos\alpha_i$ cannot be zero simultaneously, that is, $\bar{\alpha}_w$ always exists. When  $n=2$ and $\alpha_2=\pi$,  $\bar{\alpha}_w$ does not exist if $w_1=w_2$.
\begin{longlist}[(iii)]
  \item[(i)] Consider $\alpha_n<\pi/2$. Since $\tan(\cdot )$ is an increasing function in $[0,\pi/2)$ and $w_1+ \sum_{i=2}^n\cos\alpha_i$ is positive, $\bar{\alpha}_w\leq \alpha_n$ is equivalent to
\begin{eqnarray}
\label{4_convex_1} \sum_{i=2}^{n-1}w_i
\sin\alpha_i\leq \tan\alpha_n \Biggl(w_1+
\sum_{i=2}^{n-1}w_i\cos
\alpha_i \Biggr).
\end{eqnarray}
 Induction method is used to prove the  inequality (\ref{4_convex_1}). For $n=2$, (\ref{4_convex_1}) is reduced to
 \mbox{$w_1\tan\alpha_2\geq 0$} which  always holds.  We will show that  (\ref{4_convex_1}) is true for $n=k+1$ after using   it for $n=k$.  Thus,
\begin{eqnarray*}
\tan\alpha_{k+1} \Biggl(w_1+\sum
_{i=2}^{k}w_i\cos\alpha_i
\Biggr)\geq\tan\alpha_{k} \Biggl(w_1+\sum
_{i=2}^{k}w_i\cos\alpha_i
\Biggr) = \sum_{i=1}^{k}w_i\sin
\alpha_i.
\end{eqnarray*}
The above steps follow since $\tan\alpha_{k}\leq  \tan\alpha_{k+1}$ and (\ref{4_convex_1}) is assumed to be true for $n=k$.
        \item[(ii)] Now  consider the case when  $\pi/2<\alpha_n<\pi$. Since $w_1+\sum_{i=2}^{n}w_i\cos\alpha_i<0$ and $\tan(\cdot )$ is increasing  in $(\pi/2,\pi]$, $\bar{\alpha}_w\leq\alpha_n$ is equivalent to the  reverse of inequality (\ref{4_convex_1}). The proof can be completed as above.
\item[(iii)]  Cases when $\alpha_n=\pi/2$ and $\alpha_n=\pi$ are straightforward.

  \end{longlist}
Hence, the proposition is established  for $A=K_1[\aaa,\bb]$ with $\aaa,\bb\in\mathbb{T}$. In a  similar manner, it can be proved when $A$ is of types $ K_1(\aaa,\bb), K_1(\aaa,\bb]$ and $K_1[\aaa,\bb)$.

\subsection{Proof of  Proposition \texorpdfstring{\protect\ref{jen}}{2.2}}
 The statement is trivially true  when  $A=\mathbb{T}$.  As in   the proof of Proposition~\ref{convex_pr}, it is enough to prove the result only for  $A=[0,\gamma]$ with $\gamma\in(0,\pi]$. From (\ref{atan}), $\mathsf{CE}({\ttt})\geq 0$, if it exists.  Every interval is convex subset of $\mathbb{R}$.   Using the fact that an  analogous result is true for the expectation of random variables in $\mathbb{R}$   such as $\sin\ttt$ and $\cos\ttt$ (\cite{f67}, page~74),   we  can observe ranges of both $\mathsf{E}\sin\ttt$ and $\mathsf{E}\cos\ttt$.

 Suppose that  $\Pr (\ttt\in  A )=1$   and so $\Pr (\gamma-\ttt\in A )=1$. Since $\sin(\cdot )$ is nonnegative  in the range $[0,\pi]$,  $\Pr(\sin(\gamma-\ttt)\geq0)=1$. Thus,
\begin{eqnarray}
\label{e1} \mathsf{E}\sin(\gamma-\ttt)&\geq&0.
\end{eqnarray}
%

(i)  Consider the case $0<\gamma<{\pi}/{2}$. Since $\Pr(\cos\ttt>0)=1$, $\mathsf{E}\cos\ttt>0$. Dividing  both sides of (\ref{e1}) by $\cos\ttt (\mathsf{E}\cos\ttt)$ (term is positive with probability one), we get
\begin{eqnarray}
\label{zz22} ({\mathsf{E}\sin\ttt}/{\mathsf{E} \cos\ttt})\leq\tan\gamma.
\end{eqnarray}
 Using  facts that $\tan^{-1}(\cdot )$ is  increasing  in $ [0,\infty)$ and   $\mathsf{atan}(\cdot )=\tan^{-1}(\cdot )$ in the range $ [0,{\pi}/{2} )$, we  obtain $\mathsf{CE}({\ttt})\leq\gamma$.

(ii) Consider  $\gamma={\pi}/{2}$. Since both $\mathsf{E}\sin\ttt$ and $\mathsf{E}\cos\ttt$ are nonnegative and   cannot be zero simultaneously,    $\mathsf{CE}({\ttt})\leq{\pi}/{2}$.

(iii) If  ${\pi}/{2}<\gamma<\pi$, $\Pr(\cos\ttt<0)=1$ and so $\mathsf{E}\cos\ttt$ is negative. Dividing  both sides of (\ref{e1}) by $\cos\ttt (\mathsf{E}\cos\ttt)$ (positive quantity), we again obtain (\ref{zz22}). As $\tan^{-1}(\cdot )$ is  increasing in the range $(-\infty,0]$, we have
\begin{eqnarray*}
\tan^{-1} ({\mathsf{E}\sin\ttt}/{\mathsf{E}\cos\ttt} )&\leq&\gamma-
\pi.
\end{eqnarray*}
Since $\mathsf{atan}(\cdot )=\tan^{-1}(\cdot )+\pi$, $\mathsf{CE}({\ttt})\leq\gamma$.

(iv) Now consider that $\gamma=\pi$. Both $\sin\ttt\in[0,1]$ and $\cos\ttt\in[-1,1]$ with probability one. Hence $\mathsf{CE}({\ttt})\leq\pi$, if it exists.\vspace*{6pt}
%

Existence of  $\mathsf{CE}({\ttt})$ is not confirmed only for the case (iv). Suppose that  $\Pr(\ttt=0)=\Pr(\ttt=\pi)=1/2$, both $\mathsf{E}\sin\ttt$ and $\mathsf{E}\cos\ttt$ are zero, so   $\mathsf{CE}({\ttt})$  does not exist. Moreover, if $A=\mathbb{T}$ and  $\ttt$ follows uniform distribution on $\mathbb{T}$ with the density $(2\pi)^{-1}$,   $\mathsf{CE}({\ttt})$ does not exist.

\subsection{Proof of Lemma \texorpdfstring{\protect\ref{proj1}}{2.5}}
Using Lemma~\ref{rot}, it is enough to consider the case $A=[0,b]$ with $b\in(2\pi/3,\pi]$.  Decompose  the set $A$ as $A_1\cup A_2\cup A_3$, where
\begin{eqnarray*}
A_1= [0,{3b}/{4}-{\pi}/{2} ),\quad A_2=
[{3b}/{4}-{\pi}/{2},{b}/{2} ],\qquad A_3= ({b}/{2},b ],
\end{eqnarray*}
and complement of  $A$  as $B_1\cup B_2$, where
\begin{eqnarray*}
B_1= (b,\pi+{b}/{2} ],\qquad B_2= (\pi+{b}/{2},2\pi ).
\end{eqnarray*}
We have to show  for all $\psi\in A$,
\begin{eqnarray}
\label{aaa} a(\psi)=\mathsf{E}_{\psi}^{\ttt} \bigl\{u(\ttt)
\bigr\}= \int_{\ttt\in B_{1}\cup B_2}u(\ttt)f_1(\ttt\mid \psi) \,{ \rm{d}}\ttt>0,
\end{eqnarray}
where $f_1(\ttt\mid \psi)$ is the probability density  of   $\ttt$  with respect to  Lebesgue  measure and    $u(\ttt)$ is given by
\begin{eqnarray*}
u(\ttt)=\cos(\ttt_0-\psi)-\cos(\ttt-\psi),\qquad {\ttt\notin A,
\psi\in A},
\end{eqnarray*}
with    $\ttt_0$ as the projection of $\ttt$ on $A$. According to the assumption, the distribution of $\ttt$ is symmetric about $\psi$, therefore, we have
\begin{eqnarray*}
\mathsf{E}_{\psi}^{\ttt}\cos(\ttt_0-\psi)&=&
\mathsf{E}_{b-\psi}^{\ttt}\cos(\ttt_0-b+\psi).
\end{eqnarray*}
Hence, $a(\psi)=a(b-\psi)$ for $\psi\in A$. It is enough  to prove that $a(\psi)>0$  for $\psi\in A_1\cup A_2$.

Next, we examine the sign of the function  $u(\cdot )$.  Define $u(\ttt)=u_1(\ttt)+u_2(\ttt)$, where   from definition of   projection $\ttt_0$, $u_1(\ttt)$ and $u_2(\ttt)$ are given by
\begin{eqnarray*}
u_1(\ttt)&=& \cases{\displaystyle 2\sin \biggl(\frac{\ttt-b}{2} \biggr)\sin
\biggl(\frac{\ttt+b}{2}-\psi \biggr), &\quad if $\ttt \in B_1$;
\cr
0, &\quad if $\ttt \in B_2$;}
\\
u_2(\ttt)&=& \cases{ 0, &\quad if $\ttt \in B_1$;
\cr
\displaystyle 2
\sin \biggl(\frac{\ttt}{2} \biggr)\sin \biggl(\frac{\ttt}{2}-\psi \biggr),
&\quad if $\ttt \in B_2$.}
\end{eqnarray*}
If $\ttt \in B_1$,  $(\ttt-b)\in(0,{\pi}-{b}/{2}]\subset(0,2{\pi}/3)$. Hence, sign of  $u_1(\ttt)$ is only dependent on that of $\sin((\ttt+b)/2-\psi)$.  When $\psi\in A_1$, decompose  $B_1$ as  $B_1= B_{11}\cup B_{12}$, where
\begin{eqnarray*}
B_{11}=(b,2\pi+2\psi-b],\qquad B_{12}= (2\pi+2\psi-b,
\pi+{b}/{2} ].
\end{eqnarray*}
It can be noted that
\begin{eqnarray*}
\biggl(\frac{\ttt+b}{2}-\psi \biggr)\in \cases{
(b-\psi,\pi]\subset ({2
\pi}/{3},\pi ], &\quad if $\ttt \in B_{11}$;
\cr
(\pi,{\pi}/{2}+{3b}/{4}-
\psi ]\subset (\pi,{5\pi}/{4} ], &\quad if $\ttt \in B_{12}$.}
\end{eqnarray*}
Thus, when $\psi\in A_1$, $u_1(\ttt)\geq0$ for $\ttt \in B_{11}$ and  $u_1(\ttt)<0$ for $\ttt \in B_{12}$.  If  $\ttt \in B_1$ and $\psi\in A_2$, $((\ttt+b)/2-\psi)\in (\pi/3,\pi]$ and so $u_1(\ttt)\geq0$.  Similarly, when $\ttt \in B_2$, $\ttt/2\subset(2\pi/3,\pi)$,  the sign of  $u_2(\ttt)$ is only dependent on that of $\sin(\ttt/2-\psi)$. If $\ttt \in B_2$ and $\psi\in A_1\cup A_2$,  $({\ttt}/{2}-\psi)\in (\pi/4,\pi)$. This implies that  $u_2(\ttt)>0$ when $\psi\in A_1\cup A_2$.

From Table~\ref{tu1},   it is sufficient to prove
\begin{eqnarray}
\label{zz23} \int_{\ttt\in B_{12}\cup B_2}u(\ttt)f_1(\ttt\mid
\psi) \,{\rm{d}}\ttt>0\qquad\mbox{for }\psi\in A_1.
\end{eqnarray}
Now we examine two cases separately, when distribution of $\ttt$ is unimodal and bimodal.
\begin{longlist}[(ii)]
  \item[(i)] If density $f_1(\ttt\mid \psi)$ is either  $f(\ttt\mid \psi)$ or $\varepsilon f(\ttt\mid \psi)+(1-\varepsilon)f(\ttt\mid \psi+\pi)$ with   $\varepsilon\geq(1+\zeta_{\min})^{-1}$, distribution of  $\ttt$ is unimodal with mode at $\psi$. It means that  $f_1(\ttt\mid \psi)$ is increasing in   $\ttt\in[\psi+\pi,2\pi)$.    Write $u(\ttt)=v_1(\ttt)+v_2(\ttt)$, where   for $i=1,2$, $v_i(\ttt)=u_i(\ttt), \ttt\in B_{12}\cup B_2$. Thus,
\begin{eqnarray*}
\frac{v_1(\ttt)}{v_2(\ttt)}= \cases{ -\infty, &\quad if $\ttt \in B_{12}$;
\cr
0, &\quad if $\ttt \in B_2$.}
\end{eqnarray*}
Therefore,  both ${v_1(\ttt)}/{v_2(\ttt)}$ and $f_1(\ttt\mid \psi)$  are increasing in $\ttt\in B_{12}\cup B_2$. Let $\phi$ be a uniform distributed random variable  with respect to Lebesgue  measure such that $\Pr(\phi\in B_{12}\cup B_2)=1$. As $\mathsf{E}^{\phi}\{v_2(\phi)f_1(\phi\mid \psi)\}>0$,  using a result of
\cite{b84}, Theorem 2.1, we have
\begin{eqnarray*}
\frac{\mathsf{E}^{\phi}\{v_1(\phi)\}}{\mathsf{E}^{\phi}\{v_2(\phi)\}}\leq \frac{\mathsf{E}^{\phi}\{v_1(\phi)f_1(\phi\mid \psi)\}}{\mathsf{E}^{\phi}\{v_2(\phi)f_1(\phi\mid \psi)\}}.
\end{eqnarray*}
In order to prove  (\ref{zz23}), it remains to show
\begin{eqnarray}
\label{zz24} \mathsf{E}^{\phi} \bigl\{v_1(
\phi)+v_2(\phi) \bigr\}>0\qquad\mbox{for all }\psi\in
A_1.
\end{eqnarray}
The above holds since we have
\begin{eqnarray*}
(b-2\psi)\mathsf{E}^{\phi}\bigl\{v_1(\phi)+v_2(
\phi)\bigr\}&=&\int_{\ttt\in B_{12}\cup B_2} \bigl\{\cos (\phi_0-
\psi )- \cos (\phi-\psi ) \bigr\} \,{\rm{d}}\phi>0.
\end{eqnarray*}
 This completes the proof when distribution of $\ttt$ is unimodal.

\begin{table}[t]
\tabcolsep=0pt
\tablewidth=250pt
\caption{Behaviour of functions $u_1(\ttt)$ and $u_2(\ttt)$}\label{tu1}
\begin{tabular*}{\tablewidth}{@{\extracolsep{\fill}}@{}llll@{}}
\hline
& $\bolds{A_1}$ & $\bolds{A_2}$ & $\bolds{A_1\cup A_2}$\\
\hline
$B_1$ & $u_1(\ttt)\geq0$& $u_1(\ttt)\geq0$ & $u_2(\ttt)=0$\\
&$u_1(\ttt)<0$
\\[3pt]
$B_2$&$u_1(\ttt)=0$&$u_1(\ttt)=0$ & $u_2(\ttt)>0$  \\
\hline
\end{tabular*}
\end{table}

  \item[(ii)] Now we assume that $\ttt$ has a mixture distribution with the probability density  $f_1(\ttt\mid \psi)=\varepsilon f(\ttt\mid \psi)+(1-\varepsilon)f(\ttt\mid \psi+\pi)$ with $1/2\leq\varepsilon<(1+\zeta_{\min})^{-1}$, that is, the distribution of  $\ttt$ is  bimodal. It has two modes  $\psi$ and $\psi+\pi$. Antimodes are $\psi+w$ and $\psi+2\pi-w(=\beta)$, where $w=\cos^{-1} (\zeta^{-1} ({(1-\varepsilon)}/{\varepsilon} ) )$. According to condition \textup{(C2)} ,   $w\in[\pi/2,\pi)$. We may note that $\beta\leq 2\pi$ since $\psi<\pi/2\leq w$ as $\psi\in A_1$. This implies that $f_1(\ttt\mid \psi)$ is increasing in   $\ttt\in[\beta,2\pi)$. If $\beta\leq2\pi+2\psi-b$,  (\ref{zz23}) can be proved following the lines of the above case. When $\beta>2\pi+2\psi-b$, decompose the set $B_{12}\cup B_2=(2\pi+2\psi-b,2\pi)$ as $C_1\cup C_2\cup C_3$, where
\begin{eqnarray*}
C_1=(2\pi+2\psi-b,\beta],\qquad C_2=(\beta,2\pi-2w+b],
\qquad C_3=(2\pi-2w+b,2\pi).
\end{eqnarray*}
Since the probability density  $f_1(\ttt\mid \psi)$ is symmetric about $\beta$, that is, $f_1(\ttt\mid \psi)=f_1(2\beta-\ttt\mid \psi)$, we have
\begin{eqnarray*}
\int_{\ttt\in C_1}u(\ttt)f_1(\ttt\mid \psi) \,{\rm{d}}
\ttt=\int_{\ttt\in C_2} u(2\beta-\ttt)f_1(2\beta-\ttt\mid
\psi) \,{\rm{d}}\ttt=\int_{\ttt\in C_2} u(2\beta-\ttt)f_1(
\ttt\mid \psi) \,{\rm{d}}\ttt.
\end{eqnarray*}
Define
\begin{eqnarray*}
{u_3(\ttt)}= \cases{ u(\ttt)+u(2\beta-\ttt), &\quad if $\ttt \in
C_{2}$;
\cr
u(\ttt), &\quad if $\ttt \in C_3$.}
\end{eqnarray*}
 It may be noted that (\ref{zz24}) yields the following for $\psi\in A_1$:
\begin{eqnarray}
\label{bb1} \int_{\ttt\in C_{2}\cup C_3}u_3(\ttt) \,{\rm{d}}
\ttt= \int_{\ttt\in B_{12}\cup B_2}u(\ttt) \,{\rm{d}}\ttt\geq0.
\end{eqnarray}
Consider the function
\begin{eqnarray*}
u(\ttt)+u(2\beta-\ttt)=4\sin \biggl(\frac{\ttt-\ttt_0}{2} \biggr)\cos \biggl(
\frac{\ttt+\ttt_0}{2}-\beta \biggr)\sin (\beta-\psi ).
\end{eqnarray*}
When $\ttt\in B_{12}\cup B_2$, $\sin ({(\ttt-\ttt_0)}/{2} )>0$. Note that $\sin (\beta-\psi )=-\sin (w )<0$ as   $w\in[\pi/2,\pi)$. The sign of $u(\ttt)+u(2\beta-\ttt)$ is the opposite of that of $\cos ({(\ttt+\ttt_0)}/{2}-\beta )$. There are three cases according to $\pi+b/2\in C_i$, for $i=1,2,3$.  As in the case when distribution of $\ttt$ is unimodal, we define two functions $v_1(\ttt)$ and $v_2(\ttt)$ for $\ttt \in C_{2}\cup C_{3}$ in all these three cases such that  ${v_1(\ttt)}/{v_2(\ttt)}$ is increasing in $\ttt \in C_{2}\cup C_{3}$ and $v_2(\ttt)$ is nonnegative for all $\ttt \in C_{2}\cup C_{3}$.

When $\pi+b/2\in C_1$, choices are
\begin{eqnarray*}
{v_1(\ttt)}= \cases{ 0, &\quad if $\ttt \in C_{2}$;
\cr
u(\ttt), &\quad if $\ttt \in C_3$;} \qquad {v_2(\ttt)}=
\cases{ u(\ttt)+u(2\beta-\ttt), &\quad if $\ttt \in C_{2}$;
\cr
0, &
\quad if $\ttt \in C_3$.}
\end{eqnarray*}
Now consider $\pi+b/2\in C_2$. Decompose the interval  $C_2=C_{21}\cup C_{22}$, where
\begin{eqnarray*}
C_{21}=(\beta,\pi+b/2],\qquad C_{22}=(\pi+b/2,2\pi-2w+b].
\end{eqnarray*}
In this case, we choose
\begin{eqnarray*}
{v_1(\ttt)}&=& \cases{ u(\ttt)+u(2\beta-\ttt), &\quad if $\ttt \in
C_{21}$;
\cr
0, &\quad if $\ttt \in C_{22}$;
\cr
0, &\quad if
$\ttt \in C_3$;}
\\
{v_2(\ttt)}
&=& \cases{ 0, &\quad
if $\ttt \in C_{21}$;
\cr
u(\ttt)+u(2\beta-\ttt), &\quad if $\ttt \in
C_{22}$;
\cr
u(\ttt), &\quad if $\ttt \in C_3$.}
\end{eqnarray*}
When  $\pi+b/2\in C_3$, decompose the interval  $C_3=C_{31}\cup C_{32}$, where
\begin{eqnarray*}
C_{31}=(2\pi-2w+b,\pi+b/2],\qquad C_{32}=(\pi+b/2,2\pi).
\end{eqnarray*}
In this case, we define
\begin{eqnarray*}
{v_1(\ttt)}= \cases{ u(\ttt)+u(2\beta-\ttt), &\quad if $\ttt \in
C_{2}$;
\cr
0, &\quad if $\ttt \in C_{31}$;
\cr
u(\ttt), &
\quad if $\ttt \in C_{32}$;} \qquad {v_2(\ttt)}= \cases{ 0,
&\quad if $\ttt \in C_{2}$;
\cr
-u(\ttt), &\quad if $\ttt \in
C_{31}$;
\cr
0, &\quad if $\ttt \in C_{32}$.}
\end{eqnarray*}
Since  the density $f_1(\ttt\mid \psi)$ is   increasing in $\ttt \in C_{2}\cup C_{3}$,    Theorem 2.1 of \cite{b84}  completes the proof when distribution of $\ttt$ is bimodal.
\end{longlist}

\subsection{Proof of Lemma \texorpdfstring{\protect\ref{lem}}{4.1}}

See \cite{moors81}, Lemma 2, for the proof of~(i). Note that the  proof  given in \cite{moors81} utilizes the measure perseverance of the element $\mathbf{g}\in{\mathcal{G}}$ under the measure $\eta$, that is,   $\eta(\mathbf{g}^{-1}(B))\neq\eta(B)$ for all $B\in\mathfrak{B}(\mathfrak{Z})$, where  $\mathfrak{B}(\mathfrak{Z})$ consists of Borel sets  of the sample space $\mathfrak{Z}$.  Later,
Moors and van Houwelingen
\cite{moors93} relaxed this  condition of measure perseverance.

For any $\mathbf{g}\in\mathcal{G}$, if $\delta$ is an $\mathcal{G}$-equivariant estimator, we have $\mathsf{L} (\bar{\mathbf{g}}(\bolds{\nu}), \tilde{g}(\delta) )=\mathsf{L} (\bolds{\nu}, {\delta} )$, that is,
\begin{eqnarray}
\cos \bigl(h\bar{\mathbf{g}}(\bolds{\nu})-\tilde{g}({\delta}) \bigr)&=& \cos
\bigl({h}(\bolds{\nu})-{\delta} \bigr)\label{lem2},
\end{eqnarray}
for all $\bolds{\nu}\in\Omega$ and ${\delta}\in\mathcal{A}$. Substituting  ${\delta}={h}(\bolds{\nu})$, we obtain $\cos(h\bar{\mathbf{g}}(\bolds{\nu})-\tilde{g}h(\bolds{\nu}))=1$. Thus, we have $h\bar{\mathbf{g}}(\bolds{\nu})=\tilde{g}h(\bolds{\nu})\bmod(2\pi)$, for all $\nu\in\Omega$.  This proves~(ii). Exploiting  this result, (\ref{lem2}) reduces to
\begin{eqnarray*}
d \bigl(\tilde{g} {h}(\bolds{\nu}), \tilde{g}({\delta}) \bigr)&=&d \bigl({h}(
\bolds{\nu}), {\delta} \bigr),
\end{eqnarray*}
for all $\bolds{\nu}\in\Omega$ and ${\delta}\in {h}(\Omega)$. This implies that  $\tilde{{g}}(\cdot )$  is distance-preserving map on $\mathcal{A}$. Let ${\psi}$ be the projection of $\tilde{g}({\phi})$ on $\tilde{g}(A)$. Since  $\tilde{{g}}(\cdot )$ is injective, we have
\begin{eqnarray*}
d\bigl(\tilde{g}({\phi}), {\psi}\bigr)=d(\phi,\phi_{0})=d\bigl(
\tilde{g}(\phi),\tilde{g}(\phi_0)\bigr).
\end{eqnarray*}
From the uniqueness of projection, $\psi=\tilde{g}(\phi_0)$. This proves~(iii).

\subsection{Proof of Lemma \texorpdfstring{\protect\ref{ac}}{4.1}}

First, we  show that for all $\mathbf{g} \in \mathcal{G}$ and $\bolds{\nu}\in\Omega$, new estimand   ${h}_{\mathbf{z}}(\bolds{\nu})$  satisfies
\begin{eqnarray}
\label{yy} {h}_{\mathbf{g}(\mathbf{z})}\bar{\mathbf{g}}(\bolds{\nu}) = \tilde{g}
{h}_{\mathbf{z}}(\bolds{\nu})\qquad\mbox{for all }\bolds{\nu} \in \Omega.
\end{eqnarray}
Note that
\begin{eqnarray*}
\int_{\mathcal{G}} f \bigl(\mathbf{g}(\mathbf{z})\mid \bar{
\mathbf{g}}^*\bar{\mathbf{g}}(\bolds{\nu}) \bigr) \,\mathrm{d}\lambda\bigl(
\mathbf{g}^*\bigr)&=& \int_{\mathcal{G}} f \bigl(\mathbf{z}\mid \bar{
\mathbf{g}}^{-1}\bar{\mathbf{g}}^*\bar{\mathbf{g}}(\bolds{\nu}) \bigr)
\,\mathrm{d}\lambda\bigl(\mathbf{g}^*\bigr) \qquad\bigl(\mbox{from Lemma~\ref{lem}(i)}\bigr)
\\
&=&\int_{\mathcal{G}} f \bigl(\mathbf{z}\mid \bar{\mathbf{g}}^*(
\bolds{\nu}) \bigr) \,\mathrm{d}\lambda\bigl(\mathbf{g}^*\bigr)\qquad\bigl(\mbox{using transformation } \mathbf{g}^*\rightarrow\mathbf{g}\mathbf{g}^*
\mathbf{g}^{-1}\bigr).
\end{eqnarray*}
This implies that for  $\int_{\mathcal{G}} f (\mathbf{z}\mid \bar{\mathbf{g}}^*(\bolds{\nu}) ) \,\mathrm{d}\lambda(\mathbf{g}^*)=0$, (\ref{yy}) follows from Lemma~\ref{lem}(ii).  In the case of  $\int_{\mathcal{G}} f (\mathbf{z}\mid \bar{\mathbf{g}}^*(\bolds{\nu}) ) \,\mathrm{d}\lambda(\mathbf{g}^*)>0$, for any $\mathbf{g} \in \mathcal{G}$, we have
\begin{eqnarray*}
{h}_{\mathbf{g}(\mathbf{z})}\bar{\mathbf{g}} (\bolds{\nu})&=& \mathsf{atan} \biggl(
\frac{ \int_{\mathcal{G}} \sin  (\tilde{g}^*{h}\bar{\mathbf{g}}(\bolds{\nu}) ) \tau (\mathbf{g}(\mathbf{z})\mid \bar{\mathbf{g}}^*\bar{\mathbf{g}}(\bolds{\nu}) )  \,\mathrm{d}\lambda(\mathbf{g}^*)}{\int_{\mathcal{G}} \cos  (\tilde{g}^*{h}\bar{\mathbf{g}}(\bolds{\nu})  ) \tau (\mathbf{g}(\mathbf{z})\mid \bar{\mathbf{g}}^*\bar{\mathbf{g}}(\bolds{\nu}) )  \,\mathrm{d}\lambda(\mathbf{g}^*)} \biggr)
\\
&=& \mathsf{atan} \biggl(\frac{ \int_{\mathcal{G}} \sin  ({h}\bar{\mathbf{g}}^*\bar{\mathbf{g}}(\bolds{\nu})  ) \tau (\mathbf{z}\mid \bar{\mathbf{g}}^{-1}\bar{\mathbf{g}}^*\bar{\mathbf{g}}(\bolds{\nu}) )  \,\mathrm{d}\lambda(\mathbf{g}^*)}{\int_{\mathcal{G}} \cos  ({h}\bar{\mathbf{g}}^*\bar{\mathbf{g}}(\bolds{\nu})  ) \tau (\mathbf{z}\mid \bar{\mathbf{g}}^{-1}\bar{\mathbf{g}}^*\bar{\mathbf{g}}(\bolds{\nu}) )  \,\mathrm{d}\lambda(\mathbf{g}^*)} \biggr)
\\
&=& \mathsf{atan} \biggl(\frac{ \int_{\mathcal{G}} \sin  (\tilde{g}\tilde{g}^*{h}(\bolds{\nu}) ) \tau (\mathbf{z}\mid \bar{\mathbf{g}}^*(\bolds{\nu}) )  \,\mathrm{d}\lambda(\mathbf{g}^*)}{ \int_{\mathcal{G}} \cos (\tilde{g}\tilde{g}^*{h}(\bolds{\nu}) ) \tau (\mathbf{z}\mid \bar{\mathbf{g}}^*(\bolds{\nu}) )  \,\mathrm{d}\lambda(\mathbf{g}^*)} \biggr)
\\
&=& \tilde{g} \biggl(\mathsf{atan} \biggl(\frac{ \int_{\mathcal{G}} \sin  (\tilde{g}^*{h}(\bolds{\nu}) ) \tau (\mathbf{z}\mid \bar{\mathbf{g}}^*(\bolds{\nu}) )  \,\mathrm{d}\lambda(\mathbf{g}^*)}{ \int_{\mathcal{G}} \cos (\tilde{g}^*{h}(\bolds{\nu}) ) \tau (\mathbf{z}\mid \bar{\mathbf{g}}^*(\bolds{\nu}) )  \,\mathrm{d}\lambda(\mathbf{g}^*)} \biggr) \biggr)
\\
&=&\tilde{g} {h}_{\mathbf{z}}(\bolds{\nu})
\end{eqnarray*}
for all  $\bolds{\nu}\in \Omega$. The above equalities  utilize  Lemmas~\ref{lem} and~\ref{propp1},  the transformation $\mathbf{g}^*\rightarrow\mathbf{g}\mathbf{g}^*\mathbf{g}^{-1}$ and   the circular property of $\tilde{g}$. Therefore, the surjection property of  $\bar{\mathbf{g}}$ and (\ref{yy}) imply that
\begin{eqnarray}
\label{yy12} \tilde{g} {h}_{\mathbf{z}}(\Omega)= {h}_{\mathbf{g}(\mathbf{z})}(
\Omega),
\end{eqnarray}
 or equivalently, $\mathsf{cc} (\tilde{{g}}{h}_{\mathbf{z}}(\Omega) )= \mathsf{cc} ( {h}_{\mathbf{g}(\mathbf{z})}(\Omega) ) = \mathcal{A}_{\mathbf{g}(\mathbf{z})}$. Clearly, $\tilde{{g}}$-image of a closed convex set is again a closed convex set.  Therefore, $ \tilde{{g}} (\mathsf{cc} ({h}_{\mathbf{z}}(\Omega) ) )$ is also a closed convex and $\tilde{{g}} {h}_{\mathbf{z}}(\Omega)\subset \tilde{{g}} (\mathsf{cc} ({h}_{\mathbf{z}}(\Omega) ) )$. This implies that
\begin{eqnarray*}
\mathcal{A}_{\mathbf{g}(\mathbf{z})}=\mathsf{cc} \bigl(\tilde{{g}} {h}_{\mathbf{z}}(
\Omega) \bigr)\subset\tilde{{g}} \bigl(\mathsf{cc} \bigl({h}_{\mathbf{z}}(
\Omega) \bigr) \bigr)=\tilde{{g}}(A_{\mathbf{z}})
\end{eqnarray*}
 as  $\mathsf{cc} (\tilde{{g}} {h}_{\mathbf{z}}(\Omega) )$ is the  smallest convex set containing $\tilde{{g}} {h}_{\mathbf{z}}(\Omega)$.  Next, we show that $\tilde{{g}} (A_{\mathbf{z}})\subset\mathcal{A}_{\mathbf{g}(\mathbf{z})}$. As $\phi\in A_{\mathbf{z}}=\mathsf{conv}  ({h}_{\mathbf{z}}(\Omega) )\cup\mbox{bd} ( {h}_{\mathbf{z}}(\Omega) )$, there can be the following two cases:
\begin{longlist}[(ii)]
  \item[(i)] Suppose $\phi\in \mathsf{conv}  ({h}_{\mathbf{z}}(\Omega) )$. From Proposition~\ref{convex_4}, there exists $\phi_1,\phi_2\in{h}_{\mathbf{z}}(\Omega)$ such that
\begin{eqnarray*}
\phi=\mathsf{atan} \biggl(\frac{w\sin\phi_1+(1-w)\sin\phi_2}{w\cos\phi_1+(1-w)\cos\phi_2} \biggr).
\end{eqnarray*}
Operating $\tilde{g}$ on the both sides of the above equation and using the circular property, we get
\begin{eqnarray*}
\tilde{g}(\phi)=\mathsf{atan} \biggl(\frac{w\sin\tilde{g}(\phi_1)+(1-w)\sin\tilde{g}(\phi_2)}{w\cos\tilde{g}(\phi_1)+(1-w)\cos\tilde{g}(\phi_2)} \biggr).
\end{eqnarray*}
Note that $\tilde{g}(\phi_1),\tilde{g}(\phi_2) \in \tilde{g}{h}_{\mathbf{z}}(\Omega)$. Using  (\ref{yy12}), both belong to ${h}_{\mathbf{g}(\mathbf{z})}(\Omega)$. Therefore,  $\tilde{g}(\phi)\in \mathsf{conv} ({h}_{\mathbf{g}(\mathbf{z})}(\Omega) )\subset\mathcal{A}_{\mathbf{g}(\mathbf{z})}$.
  \item[(ii)] Any $\phi\in \mbox{bd} ({h}_{\mathbf{z}}(\Omega) )$ is the limit point of a series of points $\{\phi_n\}$ with $\phi_n\in \mathsf{conv} ( {h}_{\mathbf{z}}(\Omega) )$. Since $\tilde{g}$ is circular so is continuous, $\tilde{g}(\phi)$ is the limit point of the series $\{\tilde{g}(\phi_n)\}$ in  $\mathsf{conv} ( \tilde{g}{h}_{\mathbf{z}}(\Omega) )=\mathsf{conv} ({h}_{\mathbf{g}(\mathbf{z})}(\Omega) )$. Hence, $\tilde{g}(\phi)\in \operatorname{bd} ({h}_{\mathbf{g}(\mathbf{z})}(\Omega) )\subset\mathcal{A}_{\mathbf{g}(\mathbf{z})}$.
\end{longlist}

\subsection{Proof of Theorem \texorpdfstring{\protect\ref{main}}{4.1}}

 Since risk of an equivariant estimator ${\delta}(\mathbf{Z})$ is constant on the orbits of $\bolds{\nu}$ (\cite{f67}, page~149), risk of the $\mathcal{G}$-equivariant estimator  ${\delta}(\mathbf{Z})$ under the loss function $\mathsf{L}$ satisfies $\mathsf{R}(\bolds{\nu},{\delta})= \int_{\mathcal{G}}  \mathsf{R}(\bar{\mathbf{g}}(\bolds{\nu}),{\delta})  \,\mathrm{d}\lambda(\mathbf{g})$.  Using this,  the risk  of  ${\delta}(\mathbf{Z})$  is given by
 \begin{eqnarray*}
\mathsf{R}(\bolds{\nu},{\delta}) &=& \int_{\mathcal{G}} \displaystyle
\int_{\mathfrak{Z}} \bigl\{1-\cos \bigl(\delta(\mathbf{z})-h\bar{
\mathbf{g}}(\bolds{\nu}) \bigr) \bigr\} f\bigl(\mathbf{z}\mid \bar{\mathbf{g}}(\bolds{
\nu})\bigr) \,\mathrm{d}\eta(\mathbf{z})) \,\mathrm{d}\lambda(\mathbf{g})
\\
&=& \displaystyle\int_{\mathfrak{Z}}\int_{\mathcal{G}}
\bigl\{1-\cos \bigl(\delta(\mathbf{z})-\tilde{g}h(\bolds{\nu}) \bigr) \bigr\} f\bigl(
\mathbf{z}\mid \bar{\mathbf{g}}(\bolds{\nu})\bigr) \,\mathrm{d}\lambda(\mathbf{g})\,
\mathrm{d}\eta(\mathbf{z}).
\end{eqnarray*}
In the above step, we utilize the interchange in order of integration and  Lemma~\ref{lem}(ii).  Note that ${\delta}_{0}(\mathbf{g}(\mathbf{z}))$ is the projection of ${\delta}(\mathbf{g}(\mathbf{z}))$ on $\mathcal{A}_{\mathbf{g}(\mathbf{z})}$, that is,  ${\delta}_{0}(\mathbf{g}(\mathbf{z}))$ is the projection of $\tilde{g}({\delta}(\mathbf{z}))$ on $\tilde{g} (\mathcal{A}_{\mathbf{z}} )$ from  invariance of ${\delta}(\mathbf{z})$ and   Lemma~\ref{ac}.  From Lemma~\ref{lem}(iii), ${\delta}_{0}(\mathbf{g}(\mathbf{z}))=\tilde{g}({\delta}_0(\mathbf{z}))$, or equivalently, ${\delta}_0(\mathbf{z})$ is also  $\mathcal{G}$-equivariant estimator. Therefore, the  above risk expression is also valid for ${\delta}_0(\mathbf{z})$.  The difference $ \mathsf{R}(\bolds{\nu},\delta)- \mathsf{R}(\bolds{\nu},\delta_0)$ is given by
 \begin{eqnarray*}
u&=&\int_{\delta(\mathbf{z})\notin\mathcal{A}_{\mathbf{z}} } \int_{\mathcal{G}} \bigl\{\cos
\bigl(\delta_0(\mathbf{z})-\tilde{g}h(\bolds{\nu}) \bigr)-\cos \bigl(
\delta(\mathbf{z})-\tilde{g}h(\bolds{\nu}) \bigr) \bigr\} f\bigl(\mathbf{z}\mid \bar{
\mathbf{g}}(\bolds{\nu})\bigr) \,\mathrm{d}\lambda(\mathbf{g})\,\mathrm{d}\eta(
\mathbf{z})
\\
&=&\int_{\delta(\mathbf{z})\notin\mathcal{A}_{\mathbf{z}} } \biggl[\mathsf{E}^{{\mathbf{g}}} \bigl\{\cos
\bigl(\delta_0(\mathbf{z})-\tilde{g}h(\bolds{\nu}) \bigr)-\cos \bigl(
\delta(\mathbf{z})-\tilde{g}h(\bolds{\nu}) \bigr) \bigr\} \int
_{\mathcal{G}} f\bigl(\mathbf{z}\mid \bar{\mathbf{g}}(\bolds{\nu})\bigr) \,
\mathrm{d}\lambda(\mathbf{g}) \biggr]\,\mathrm{d}\eta(\mathbf{z}),
\end{eqnarray*}
if $\int_{\mathcal{G}} f(\mathbf{z}\mid \bar{\mathbf{g}}(\bolds{\nu}))  \,\mathrm{d}\lambda(\mathbf{g})>0$, where the expectation is taken over $\mathbf{g}$ with respect to a probability measure $\tau (\mathbf{z}\mid \bar{\mathbf{g}}(\bolds{\nu}) )  \,\mathrm{d}\lambda(\mathbf{g})$. Using\vspace*{2pt} the representation of $h_{\mathbf{z}}(\bolds{\nu})$ given in (\ref{res_h}) and denoting by $v_{\mathbf{z}}(\bolds{\nu})= [\{\mathsf{E} \sin  (\tilde{g}{h}(\bolds{\nu}) )\}^2+\{\mathsf{E} \sin  (\tilde{g}{h}(\bolds{\nu}) )\}^2 ]^{{1}/{2}}$, we have
\begin{eqnarray*}
\sin h_{\mathbf{z}}(\bolds{\nu})&=& {\mathsf{E} \sin \bigl(\tilde{g} {h}(
\bolds{\nu}) \bigr)}/{v_{\mathbf{z}}(\bolds{\nu})},
\\
\cos h_{\mathbf{z}}(\bolds{\nu})&=&{\mathsf{E}\cos \bigl(\tilde{g} {h}(\bolds{
\nu}) \bigr)}/{v_{\mathbf{z}}(\bolds{\nu})},
\end{eqnarray*}
the risk difference  is given by
 \begin{eqnarray*}
u&=&\int_{\delta(\mathbf{z})\notin\mathcal{A}_{\mathbf{z}} } {v_{\mathbf{z}}(\bolds{\nu})} \bigl\{\cos
\bigl(\delta_{0}(\mathbf{z})-h_{\mathbf{z}}(\bolds{\nu})\bigr)-\cos
\bigl(\delta(\mathbf{z})-h_{\mathbf{z}}(\bolds{\nu})\bigr) \bigr\} \int
_{\mathcal{G}} f\bigl(\mathbf{z}\mid \bar{\mathbf{g}}(\bolds{\nu})\bigr) \,
\mathrm{d}\lambda(\mathbf{g})\,\mathrm{d}\eta(\mathbf{z}).
\end{eqnarray*}
If $l(\mathcal{A}_{\mathbf{z}})\leq (2/3)\pi$,    the above integrand is positive since  $\cos(\delta(\mathbf{z})-h_{\mathbf{z}}(\bolds{\nu}))\leq\cos(\delta_0(\mathbf{z})-h_{\mathbf{z}}(\bolds{\nu}))$ from Lemma~\ref{aaaa}.  This completes the proof.
\end{appendix}

\section*{Acknowledgements}
The authors are thankful to the two  referees and an associate editor for their constructive  comments and suggestions which have substantially improved the paper.



\printhistory
\end{document}